# Reliability Analysis of a Multi-State Truly-Threshold System Using a Multi-Valued Karnaugh Map

**Ali Muhammad Ali Rushdi and Fares Ahmad Muhammad Ghaleb**
Department of Electrical and Computer Engineering
Faculty of Engineering, King Abdulaziz University
P. O. Box 80200, Jeddah, 21589, Saudi Arabia
E-mail: arushdi@kau.edu.sa

**Abstract**
This paper is devoted to the Boolean-based analysis of a prominent class of non-repairable coherent multi-state systems with independent non-identical multi-state components. This class of systems is represented by a multi-state coherent truly-threshold system of $(m + 1)$ states, which is not necessarily binary-imaged. The paper represents such a system *via* Boolean expressions of system success or system failure at each non-zero level, which are in the form of (a) minimal sum-of-products formulas, or (b) disjoint sum-of-products formulas, which are probability–ready expressions that are directly convertible to probabilistically expected values. Several map representations are also offered, including a single multi-value Karnaugh map (MVKM) depicting system success and $M$ maps of binary contents (albeit with multi-valued input variables) depicting system success or system failure at each non-zero level. Our results demonstrate that, similarly to binary threshold systems, multi-state threshold systems possess the shellability property. In fact, the minimal sum-of-products formulas for system success or system failure at each non-zero level can be converted to disjoint formulas without any change in the number of products. This work stresses the adaptation of techniques initially addressing binary systems, for application to multi-state ones,rather than innovating novel concepts for these latter systems.

**Keywords-** System reliability, Threshold system, Multi-state system, Multi-valued Karnaugh map (MVKM), System levels, Minimal upper vector, Maximal lower vector, Shellability.

## 1. Introduction

The study of system reliability in engineering practice has naturally been confined in its beginnings to the binary model, in which both the system and its components have only two possible states: good (working) or failed (non-working). Such a confinement was necessitated by mathematical simplicity and was made possible through the dichotomization of processes that usually involved more than two states. During the past four decades, many researchers advocated a paradigm shift in system reliability studies from the model of a too-simple binary system to that of a more appealing multi-state system (MSS) [1-22], but most progress was still made solely in the binary domain, where many new reliability models continued to be introduced. A prominent such model is the threshold model [23-32], which is a system that is successful if and only if the weighted sum of its component successes equals or exceeds a certain threshold. The binary threshold model is characterized by a success/failure function that is shellable (with minimal sum-of-products formulas that can be converted to disjoint ones without any increase in the number of products) [24, 26-28, 31]. In the symmetric case (when component weights are the same), the threshold system reduces to the $k$-out-of-$n$ system $(1 \le k \le n)$, with $k$ replacing the ceiling of the threshold divided by the common weight [23]. The threshold system is commonly renamed as the weighted k-out-of-n system [32-39], wherein the symbol $k$ assumes the role of the number of components while the symbol $n$ deviates from its customary meaning to usurp the role of the threshold. With this unfortunate renaming, you can hear of the ridiculous name of a weighted 6-out-of-5 system, which defies the streamlined logic of the English language. The utility of the binary threshold model is not confined to the reliability domain, since it is very useful also in the study of weighted voting systems [40-44].

In a recent paper [45], the present authors made a serious attempt to unify the MSS and threshold models for the case of a coherent system with a well-defined binary image [46]. This paper introduces a more versatile multi-state coherent threshold model that is more likely to inherit all strengths of both the multi-state model and the threshold model. The proposed model is truly threshold in the wider multi-state sense and might prove to possess applications that are far more versatile than those of its two constituent sub-models. The new model does not require the pertinent system to be binary-imaged. In fact, it lacks a binary image for any of its levels. Though this novel system might prove somewhat more difficult to analyze, it will still yields to most techniques in our arsenal of tools for MSS systems [45].

This paper is a part of an on-going activity that strives to provide a pedagogical Boolean-based treatment of binary reliability problems [23-32, 47-55], as well as multi-state reliability ones [56-68]. We aspire to establish a clear and insightful interrelationship between the two-state modeling and the multi-state one by stressing that multi-valued concepts are natural and simple extensions of two-valued ones. Moreover, we extend algebraic techniques and tools of switching algebra or binary logic to ones of multiple-valued logic, so as to evaluate each of the multiple levels of the system output as an individual binary or propositional function of the system multi-valued inputs. The formula of each of these levels is then written in a disjoint sum-of-products form, which is a probability–ready expression, thereby allowing its immediate conversion, on a one-to-one basis, into a probability or expected value. It is well-known that visual insight secured through the use of Karnaugh maps aids dramatically in the comprehension of coherent-system concepts, whether they are binary or multi state. Therefore, the present paper strives to provide such visual insight through the use of a single multi-value Karnaugh map (MVKM) [45, 56, 57, 58, 67-75] depicting system success and M maps of binary contents (albeit with multi-valued input variables) depicting the binary system success or system failure at each non-zero level.

The organization of the remainder of this paper is as follows. Section 2 presents important assumptions, notation and specifies the general characteristics of a coherent truly-threshold multi-state system. Section 3 introduces the running example used in this paper. Section 4 develops a single multi-value Karnaugh map (MVKM) that represents system success, three maps of binary contents (albeit with multi-valued variables) depicting the binary system success at each non-zero level, and three more maps also of binary contents (but again with multi-valued variables) depicting the binary system failure at each non-zero level. These maps are used in the development of Boolean expressions of system success or system failure at each non-zero level. These expressions are in the form of (a) minimal sum-of-products formulas, or (b) disjoint sum-of-products formulas, which are probability–ready expressions that are directly convertible to expected values.  Section 5 concludes the paper. To make the paper self-contained, it is supplemented with two appendices. Appendix A deals with the extension of the concept of a probability-ready expression (PRE) from the binary to the multi-state case. Likewise, Appendix B reviews the concept of the Boole-Shannon expansion in a multi-valued sense. Work in this paper is a sequel of our earlier work in references [45, 59, 68], and we are relying herein on the detailed nomenclature reported in these references.

## 2. Assumptions, Notation and Specification of a Running Example

### 2.1 Assumptions

- The model considered is one of a multi-state system with multistate components, specified by the structure or success function $S(X)$ [13, 45, 68]

$$S: \{0, 1, \cdots, m_1\} \times \{0, 1, \cdots, m_2\} \times \ldots \times \{0, 1, \cdots, m_n\} \to \{0, 1, \cdots, M\}. \tag{1}$$

- The system is generally non-homogeneous, i.e., the number of system states $(M + 1)$ and the numbers of component states $(m_1 + 1), (m_2 + 1), \cdots, (m_n + 1)$ might differ. When these numbers have a common value, the system reduces to a homogeneous one.
- The system is a non-repairable one with statistically independent non-identical (heterogeneous) components.
- The system is a coherent one enjoying the properties of causality, monotonicity, and component relevancy [1, 2, 4, 24, 28, 45, 56, 57, 59, 67, 68].
- The system is not necessarily binary-imaged [14, 45, 46, 58, 59, 67], i.e., system success at each specific level is not necessarily dependent only on component successes at the same level.

### 2.2 Notation

| Symbol | Description |
|---|---|
| $X_k$ | A multivalued input variable representing component $k$ $(1 \leq k \leq n)$, where $X_k \in \{0, 1, \ldots, m_k\}$, and $m_k \geq 1$ is the highest value of $X_k$. |
| $X_k\{j\}$ | A binary variable representing instant $j$ of $X_k$ $$X_k\{j\} = \{X_k = j\},$$ i.e., $X_k\{j\} = 1$ if $X_k = j$ and $X_k\{j\} = 0$ if $X_k \neq j$. The instances $X_k\{j\}$ for $\{0 \leq j \leq m_k\}$ form an orthonormal set, namely, for $\{1 \leq k \leq n\}$ $$\vee_{j=0}^{m_k} X_k\{j\} = 1, \tag{2a}$$ $$X_k(j_1)\, X_k(j_2) = 0 \text{ for } j_1 \neq j_2. \tag{2b}$$ Orthonormality is very useful in constructing inverses or complements. The complement of the union of certain instances is the union of the complementary instances. In particular, the complement of $X_k\{\geq j\} = X_k\{j, j+1, \ldots, m_k\}$ is $X_k\{< j\} = X_k\{0, 1, \ldots, j-1\}$. |
| $X_k\{\geq j\}$ | An upper value of $X_k$ $\{0 \leq j \leq m_k\}$: $$X_k\{\geq j\} = X_k\{j, j+1, \ldots, m_k\} = \vee_{i=j}^{m_k} X_k\{i\} = X_k\{j\} \vee X_k\{j+1\} \vee \ldots \vee X_k\{m_k\}. \tag{3}$$ The value $X_k\{\geq 0\}$ is identically 1. The set $X_k\{\geq j\}$ for $\{1 \leq j \leq m_k\}$ is neither independent nor disjoint, and hence it is difficult to be handled mathematically, but it is very convenient for translating the verbal or map/tabular description of a coherent component into a mathematical form |

| | when viewing component success at level $j$. The complement of $X_k\{\geq j\}$ is |
|---|---|
| | $$X_k\{< j\} = X\{0, 1, \ldots, j-1\} = X_k\{0\} \vee X_k\{1\} \ldots \vee X_k\{j-1\} = X_k\{\leq (j-1)\}. \quad (4)$$ |
| $X_k\{\leq j\}$ | A lower value of $X_k$ $\{0 \leq j \leq m_k\}$: $$X_k\{\leq j\} = X_k\{0, 1, \ldots, j-1, j\} = \vee_{i=0}^{j} X_k\{i\} = X_k\{0\} \vee X_k\{1\} \ldots \vee X_k\{j-1\} \vee X_k\{j\}. \quad (5)$$ The value $X_k\{\leq m_k\}$ is identically 1. The set $X_k\{\leq j\}$ for $\{0 \leq j \leq (m_k - 1)\}$ is neither independent nor disjoint, and hence it is not convenient for mathematical manipulation though it is suitable for expressing component failure at level $(j + 1)$. Instances, upper values and lower values are related by $$X_k\{j\} = X_k\{\geq j\} X_k\{< (j+1)\} = X_k\{\geq j\} \bar{X}_k\{\geq (j+1)\} = X_k\{\leq j\} X_k\{> (j-1)\}$$ $$= X_k\{\leq j\} \bar{X}_k\{\leq (j-1)\}. \quad (6)$$ |
| $S$ | A multivalued output variable representing the system, where $$S \in \{0, 1, \ldots, M\}, \quad (7)$$ and $M \geq 1$ is the highest value attained by the system. The system is called homogeneous if $M = m_1 = m_2 = \cdots = m_n$. The function $S(X)$ is usually called the system success or the structure function. It is conveniently represented by a Multi-Valued Karnaugh Map (MVKM) [45, 56, 57, 59, 67-75]. Its complement $\bar{S}(X)$ is called system failure and is also a multivalued variable of $(M + 1)$ values. The arithmetic sum $(S(X) + \bar{S}(X))$ is identically equal to $M$. |
| $S\{j\}$ | A binary variable representing instant $j$ of $S$ $$S\{j\} = \{S(X) = j\}, \quad (8)$$ i.e., $S\{j\} = 1$ if $S(X) = j$, and $S(j) = 0$ if $S(X) \neq j$. The instances $S\{j\}$ for $\{0 \leq j \leq M\}$ form an orthonormal set, i.e. $$\vee_{j=0}^{M} S\{j\} = 1, \quad (9)$$ $$S\{j_1\} S\{j_2\} = 0 \text{ for } j_1 \neq j_2, \quad (10)$$ which means that one, and only one, of the $(M + 1)$ instances of $S$ has the value 1, while the other instances are all $0's$. |
| $S\{\geq j\}$ | An upper value of $S$ (called system success at level $j$) $$S\{\geq j\} = S\{j, j+1, \ldots, M\} = \bigvee_{i=j}^{M} S\{i\}. \quad (11)$$ |
| $S\{\leq j\}$ | A lower value of $S$ (called system failure at level $j + 1$) |

$$S\{\leq j\} = S\{< (j+1)\} = S\{0, 1, \ldots, j\} = \bigvee_{i=0}^{j} S\{i\}. \tag{12}$$

Instances, upper values and lower values of $S$ are related by

$$S\{j\} = S\{\geq j\} S\{< (j+1)\} = S\{\geq j\} \bar{S}\{\geq (j+1)\} = S\{\leq j\} S\{> (j-1)\}. \tag{13}$$

### 2.3. Specification of a Coherent Truly-Threshold Multi-State System

A coherent truly-threshold multi-state system might be defined to be in state $j$ ($0 \leq j \leq M$) or above if and only if the sum of its component successes $X_i$, each weighted by a specific non-negative weight $W_i$, is greater than or equal to a specific non-negative threshold $T_j$

$$\{S\{\geq j\} = 1\} \quad \text{iff} \quad \{ \sum_{i=1}^{n} W_i X_i \geq T_j \}, \tag{14}$$

If the condition (14) is satisfied, the system success at level $j$ ($1 \leq j \leq M$) is said to be a switching threshold function [23-32]. Correspondingly, the system is said to be exactly at level $j$ ($0 \leq j \leq M$) iff it is in state $j$ or above but it is below state $j+1$.

$$\{S\{j\} = 1\} \quad \text{iff} \quad \{ T_j \leq \sum_{i=1}^{n} W_i X_i < T_{j+1} \}, \tag{15}$$

System coherence necessitates that $T_0$ be zero, and that the successive thresholds are totally ordered as

$$0 = T_0 < T_1 < T_2 < \cdots . T_j < T_{j+1} < \cdots < T_M, \tag{16}$$

With these definitions, the system is guaranteed to be in state 0 or above $\{S\{\geq 0\} = 1\}$ since every $X_i \geq 0$ and every $W_i \geq 0$, while $T_0 = 0$. Therefore, we do not need to consider $S\{\geq 0\}$ in our further exploration. We considered a component weight $W_i$ to be independent of the system level $j$, but we might enhance modelling versatility if we utilize $W_{ij}$ instead of just $W_i$, i.e. if we take a component weight to be dependent on the system level $j$.

### 3. Description of a Running Practical Example

We employ as a specific running example a plausible extension of an example given by Wood [8]. A 4-engine aircraft is modeled as a 4-state system of four 3-state components, where the component and system states are defined as shown in Tables 1 and 2, respectively. The dependence of the multi-valued system success S on the multi-valued component successes $X$ is limited to dependence on the sum $(X_1 + X_2 + X_3 + X_4)$ of component successes. The specifications given amount to restricting the general system described in Sec. 2 to $m = 2, M = 3, W_1 = W_2 = W_3 = W_4 = 1$. We further choose

$$T_0 = 0, T_1 = 2, T_2 = 4, T_3 = 6 \text{ and } T_4 = 9. \tag{17}$$

With these choices, we obtain the relation shown in Table 3 between the system level $j$ attained and the corresponding range $[T_j, T_{j+1})$ to which the weighted sum of component successes $\sum_{i=1}^{4} X_i$ belongs.

To facilitate the characterization of the system under study, we first construct in Fig. 1 a MVKM representing the sum of the multi-valued variables $X_1, X_2, X_3$ and $X_4$. Noting that the system is in state $j$, where $0 \leq j \leq 3$, iff the sum in Fig.1 lies in $[T_j, T_{j+1})$, i.e., iff

$$T_j \leq \sum_{i=1}^{4} X_i < T_{j+1}, \qquad (18)$$

we use Fig. 1 and Table 3 to construct the MVKM in Fig. 2 for the multi-valued system success $S(X)$. In the remainder of this section, we explain the characterizing features of a minimal upper vector and a maximal lower vector and explore the properties of symmetry, coherence, and dominance of the present system, and we occasionally utilize the MVKM in Fig. 2, and its versions in Figs. 3-5, to elaborate on these feature and properties.

A minimal upper vector (MUV) at level $j > 0$, denoted $\theta_{ji}$ is an upper vector $X$ for level $j$ such that $S(Y) < j, \{j = 1, 2, ..., M\}$ for any vector $Y < X$. A maximal lower vector (MLV) at level $j < M$, denoted $\sigma_{ji}$ is a lower vector $X$ for level $j$ such that $S\{Y\} > j, \{j = 1, 2, ..., M\}$ for any vector $Y > X$. Figure 3 highlights a typical MUV, while Fig. 4 demonstrates a typical MLU for our running example. In the sequel, we explain how to obtain all MUVs and all MULs for this system. There is a subtle relation between a minimal upper vector (MUV) at a certain level and a prime implicant of success (minimal path) at that level, and, likewise, there is a dual relation between a maximal lower vector (MLV) at a certain level and a prime implicant of failure (minimal cutest) at that level. A minimal path constitutes all the upper vectors extending (inclusively) from a particular MUV to the all-highest vector, while a minimal cutset comprises all the upper vectors extending (inclusively) from the all-0 vector to a particular MLV [45, 59, 68].

The system function $S(X)$ is a **totally symmetric** function. This means that it is partially symmetric in each of the six pairs of variables $X_i$ and $X_j$, where $(i,j) = (1,2), (1,3), (1,4), (2,3), (2,4)$ or $(3,4)$. It suffices to prove partial symmetry for three out of the aforementioned six pairs, say pairs $(1,2), (3,4)$ and $(1,3)$. Partial symmetry in the variables $X_i$ and $X_j$ means that $X_i$ and $X_j$ can be interchanged without a change in the value of the function $S(X)$. This is equivalent to equality of the multivalued quotients:

$$S(X)/X_i\{k\} = S(X)/X_j\{k\}, \qquad k = 0,1,2. \qquad (19)$$

or to equality of the multivalued quotients:

$$S(X)/X_i\{k_1\}X_j\{k_2\} = S(X)/X_i\{k_2\}X_j\{k_1\}, \qquad k_1 \neq k_2 \in \{0,1,2\}. \qquad (20)$$

The MVKM demonstrates partial symmetry for $(i,j) = (1,2)$ through equivalence of columns 2 and 4, columns 3 and 7, and columns 6 and 8. Likewise, the MVKM demonstrates partial symmetry for $(i,j) = (3,4)$ through equivalence of the corresponding rows. Partial symmetry for $(i,j) = (1,3)$ or $(2,4)$ is evident since the contents of row $m$ and column $m$ are the same.

**Coherence** of the present multi-state system is illustrated neatly by the MVKM. In fact, the three properties of coherence can be observed from the MVKM of Fig. 2 (reproduced in Fig. 5) as follows:

1.   **Causality** is evident from the map since the entry of the all-0 cell is 0 (indicating that S(0, 0, 0, 0) = 0 ) and the entry of the all-2 cell is 3 (indicating that S (2, 2, 2, 2) = 3 ).

**Monotonicity** with respect to component 1 can be observed by dividing the map into three submaps each comprising one third of the total map and consisting of three map columns. These three submaps are separated by two vertical $X_1$ borders transitioning from $X_1 = 0$ to $X_1 = 1$, and from $X_1 = 1$ to $X_1 = 2$, respectively. They represent the regions of $X_1 = 0$, $X_1 = 1$, and $X_1 = 2$ (corresponding to $X_1\{0\} = 1, X_1\{1\} = 1$, and $X_1\{2\} = 1$), and they stand for the multivalued quotients $S(X)/X_1\{0\}$, $S(X)/X_1\{1\}$, and $S(X)/X_1\{2\}$. Monotonicity w.r.t. components 1 is demonstrated by the relations:

$$S(X)/X_1\{2\} \geq S(X)/X_1\{1\} \geq S(X)/X_1\{0\}. \tag{21}$$

This means that the entry of any cell in the third (rightmost) submap $S(X)/X_1\{2\}$ is equal to or greater than the entry of the corresponding cell in the second (middle) submap $S(X)/X_1\{1\}$, which, in turn, is equal to or greater than the entry of the corresponding cell in the first (leftmost) submap $S(X)/X_1\{0\}$. This asserts monotonicity w.r.t. $X_1$, i.e., for fixed values of $X_2$, $X_3$, and $X_4$, the output $S$ is non-decreasing when $X_1$ increases. The MVKM in Fig. 5 can be similarly used to demonstrate monotonicity of its function $S(X)$ w.r.t. each of the remaining components 2, 3, and 4, namely:

$$S(X)/X_k\{2\} \geq S(X)/X_k\{1\} \geq S(X)/X_k\{0\}, \quad k = 2,3,4. \tag{22}$$

Otherwise, it might suffice to attribute monotonicity of the remaining components to the fact that $S(X)$ is totally symmetric, together with monotonicity w.r.t. $X_1$.

**Relevancy** of component 1 is evident from the fact that there are cases in which the inequalities in (21) are strict ones. For example, component 1 is relevant in the transition from state 0 to state 1 since

$$S(0,1,0,0) = 0 < 1 = S(1,1,0,0), \text{ where } (0,1,0,0) < (1,1,0,0). \tag{23a}$$

Likewise, component 1 is relevant in the transition from state 1 to state 2 since

$$S(1,2,0,0) = 1 < 2 = S(2,2,0,0), \text{ where } (1,2,0,0) < (2,2,0,0). \tag{23b}$$

The MVKM in Fig. 5 can be used in a similar fashion to demonstrate relevancy of each of the components 2, 3, and 4. Otherwise, it might suffice to attribute relevancy of the remaining components to the fact that $S(X)$ is totally symmetric, together with the relevancy of $X_1$.

We recall that a dominant multi-state system is a coherent multi-state system, in which $S(X) > S(Y)$ implies vector $X$ dominates vector $Y$, which means that vector $X$ must be larger than a vector that is in the same equivalence class as vector $Y$, albeit that vector $X$ may not be larger than vector $Y$ itself. This means that in a dominant system, every vector of state $j > 0$ must be larger than at least one vector of a smaller state value [14, 59]. A dominant system can be further classified as either (a) with a binary image, or (b) without a binary image. We now identify our running example as a dominant system without a binary image. For example, Fig. 5 reveals that vectors $X = (2,2,0,2)$ and $Y = (0,0,2,2)$ are such that $S(X) = 3 > 2 = S(Y)$, and $X$ dominates $Y$ (though $X$ is incomparable to $Y$, since $Y$ in the same class as $Z = (0,2,0,2)$ with $S(Y) = S(Z)$ and $X > Z$).

A multi-state system is completely characterized by of the set of minimal upper vectors (MUVs) at level $j > 0$ $(j = M, (M-1), ..., 2, 1)$, denoted by $\boldsymbol{\theta}(j)$. This is the set of vectors $X$, which are (a) upper for level $j$ $(S(X) \geq j)$, and (b) minimal among the upper vectors (such that $S(Y) < j$, for any vector $Y < X$). This set leads to a complete non-binary-imaged characterization of the system under study in terms of the binary functions $S\{\geq j\}$ $(j = M, (M-1), ..., 2, 1)$, where, $S\{\geq j\}$ depicts system success at level $j$ in a minimal sum-of-products form [45, 59, 68]. Figure 3 demonstrates that the cell/vector $X = [1, 2, 2, 2]$ (highlighted in yellow) is an MUV for level 3 since $S(X) = 3 \geq 3$, while all the cells/vectors $Y < X$ (highlighted in green) are such that $S(Y) < 3$.

An alternative complete characterization is *via* maximal lower vectors (MLVs) at level $j < M$ $(j = 0, 1, ..., (M-1))$, denoted by $\boldsymbol{\sigma}(j)$. This is the set of vectors $X$, which are (a) lower for level $j$ $(S(X) \leq j)$, and (b) maximal among the lower vectors (such that $S(Y) > j$, for any vector $Y > X$). This set leads to a complete non-binary-imaged characterization of the system under study in terms of the binary functions $S\{\leq j\} = S\{< (j+1)\}$ $(j = 0, 1, ..., (M-1))$, where, $S\{< (j+1)\} = S\{\leq j\}$ depicts system failure at level $(j+1)$, again in a minimal sum-of-products form [45, 59, 68]. Figure 4 demonstrates that the cell/vector $X = [2, 1, 1, 1]$ (highlighted in yellow) is an MLV for level 2 since $S(X) = 2 \leq 2$, while all the cells/vectors $Y > X$ (highlighted in green) are such that $S(Y) > 2$.

We recall that a binary-imaged multi-state system is a system whose success at level $j$ is a function only of component successes at the same level ($S\{\geq j\}$ is a function of $X\{\geq j\}$ only), or equivalently, it is a system whose failure at level $j$ is a function only of component failures at the same level ($S\{\leq (j-1)\}$ is a function of $X\{\leq (j-1)\}$ only) [45, 59]. For a binary-imaged system, elements of the set of MUVs $\boldsymbol{\theta}(j)$ are vectors of $j$ or 0 components only, and elements of the set of MLVs $\boldsymbol{\sigma}(j)$ are vectors of $j$ or $M$ components only [45, 59]. In the sequel, we assert that our running system lacks the binary-imaged properties for its MUVs and MLVs, and hence it is not a binary-imaged system.

## 4. Analysis of the Running Practical Example

Table 4 provides a listing of sample elements of the set of Minimal Upper Vectors (MUVs) $\boldsymbol{\theta}(k)$, $(k = 3, 2, 1)$ through the solution of the equations:

$$X_1 + X_2 + X_3 + X_4 = 2k, \qquad k = 3, 2, 1, \qquad (24)$$

for integers $X_1, X_2, X_3, X_4$ lying in $\{0, 1, 2\}$.

Dually, Table 5 provides a listing of sample elements of the set of Maximal Lower Vectors (MLVs) $\boldsymbol{\sigma}(k)$, $(k = 2, 1, 0)$ through solution of the equations:

$$X_1 + X_2 + X_3 + X_4 = 2k + 1, \qquad k = 2, 1, 0, \qquad (25)$$

for integers $X_1, X_2, X_3, X_4$ lying in $\{0, 1, 2\}$.

Two complete solutions of the present problem are given in various variants of six consecutive figures. Variants of three consecutive figures (Figs. 6-8) present a solution in a success perspective, while variants of the next three consecutive figures (Figs. 9-11) present a dual solution in a failure perspective.

Variants of Fig. 6 deal with binary success at level 3, i.e., with $S\{\geq 3\} = S\{3\}$. Figure 6a shows the binary success at level 3 with MUV cells highlighted in yellow for MUV cells of the (2, 2, 1, 1)-type, and in green for MUV cells of the (2, 2, 2, 0)-type. Figure 6b demonstrates two sample prime-implicant (PI) loops for the binary success at level 3. One PI loop extends from the yellow (2, 2, 2, 0)-type MUV cell (0,2,2,2) to the all-2 cell (2,2,2,2), and another PI loop extends from the green (2, 2, 1, 1)-type MUV cell (1,2,1,2) again to the all-2 cell. The total number of such PI loops is c(4, 1)+ c(4, 2) = 4 + 6 = 10, where c(n, k) is the combinatorial coefficient or n choose k. The binary success at level 3 is given by the minimal sum-of-products expression:

$$S\{\geq 3\} = S\{3\} = X_1\{\geq 2\} X_2\{\geq 2\} X_3\{\geq 2\} X_4\{\geq 0\} \vee X_1\{\geq 2\} X_2\{\geq 2\} X_3\{\geq 0\} X_4\{\geq 2\} \vee X_1\{\geq 2\} X_2\{\geq 0\} X_3\{\geq 2\} X_4\{\geq 2\} \vee X_1\{\geq 0\} X_2\{\geq 2\} X_3\{\geq 2\} X_4\{\geq 2\} \vee X_1\{\geq 2\} X_2\{\geq 2\} X_3\{\geq 1\} X_4\{\geq 1\} \vee X_1\{\geq 2\} X_2\{\geq 1\} X_3\{\geq 2\} X_4\{\geq 1\} \vee X_1\{\geq 1\} X_2\{\geq 2\} X_3\{\geq 2\} X_4\{\geq 1\} \vee X_1\{\geq 2\} X_2\{\geq 1\} X_3\{\geq 1\} X_4\{\geq 2\} \vee X_1\{\geq 1\} X_2\{\geq 2\} X_3\{\geq 1\} X_4\{\geq 2\} \vee X_1\{\geq 1\} X_2\{\geq 1\} X_3\{\geq 2\} X_4\{\geq 2\}. \qquad (26)$$

Figure 6c covers the binary success at level 3 with ten disjoint (non-overlapping) loops, comprising: (a) four loops (an intact 3-cell loop plus three curtailed 2-cell loops) that inherit the ones whose MUVs are (2, 2, 2, 0)-type, and (b) six loops (each a curtailed 1-cell loop containing just the MUV) that inherit the ones whose MUVs are (2, 2, 1, 1)-type. This threshold function is a shellable one, with each of the original PI loops retained as a single disjoint loop. The PRE for the binary success at level 3 is:

$$S_{PRE}\{\geq 3\} = S_{PRE}\{3\} = X_1\{2\} X_2\{2\} X_3\{2\} X_4\{0,1,2\} \vee X_1\{2\} X_2\{2\} X_3\{0,1\} X_4\{2\} \vee X_1\{2\} X_2\{0,1\} X_3\{2\} X_4\{2\} \vee X_1\{0,1\} X_2\{2\} X_3\{2\} X_4\{2\} \vee X_1\{2\} X_2\{2\} X_3\{1\} X_4\{1\} \vee X_1\{2\} X_2\{1\} X_3\{2\} X_4\{1\} \vee X_1\{1\} X_2\{2\} X_3\{2\} X_4\{1\} \vee X_1\{2\} X_2\{1\} X_3\{1\} X_4\{2\} \vee X_1\{1\} X_2\{2\} X_3\{1\} X_4\{2\} \vee X_1\{1\} X_2\{1\} X_3\{2\} X_4\{2\}. \qquad (27)$$

Figure 6d presents another coverage for the binary success at level 3 with alternative ten disjoint (non-overlapping) loops, comprising: (a) six loops (an intact 4-cell loop plus two curtailed 2-cell loops and three curtailed 1-cell loops) that inherit the ones whose MUVs are (2, 2, 1, 1)-type, and (b) four loops (each a curtailed 1-cell loop containing just the MUV) that inherit the ones whose MUVs are (2, 2, 2, 0)-type. This confirms the earlier observation that this threshold function is a shellable one, with each of the original PI loops retained as a single disjoint loop. Therefore, another PRE for the binary success at level 3 is:

$$S_{PRE}\{\geq 3\} = S_{PRE}\{3\} = X_1\{2\} X_2\{1,2\} X_3\{2\} X_4\{1,2\} \vee X_1\{2\} X_2\{1,2\} X_3\{1\} X_4\{2\} \vee X_1\{1\} X_2\{1,2\} X_3\{2\} X_4\{2\} \vee X_1\{2\} X_2\{2\} X_3\{1\} X_4\{1\} \vee X_1\{1\} X_2\{2\} X_3\{1\} X_4\{2\} \vee X_1\{1\} X_2\{2\} X_3\{2\} X_4\{1\} \vee X_1\{2\} X_2\{2\} X_3\{2\} X_4\{0\} \vee X_1\{2\} X_2\{2\} X_3\{0\} X_4\{2\} \vee X_1\{2\} X_2\{0\} X_3\{2\} X_4\{2\} \vee X_1\{0\} X_2\{2\} X_3\{2\} X_4\{2\}. \qquad (28)$$

Variants of Fig. 7 deal with binary success at level 2, i.e., with $S\{\geq 2\} = S\{2,3\}$. Figure 7a shows the binary success at level 2 with MUV cells highlighted in yellow for MUV cells of the (2, 2, 0, 0)-type, in green for MUV cells of the (2, 1, 1, 0)-type, and finally in purple for MUV cells of the (1, 1, 1, 1)-type. Figure 7b demonstrates three sample prime-implicant (PI) loops for the binary success at level 2. These are (a) a 9-cell red PI loop extending from the yellow (2, 2, 0, 0)-type MUV cell (2,0,2,0) to the all-2 cell (2, 2, 2, 2), (b) a shaded 12-cell PI loop extending from the (2, 1, 1, 0)-type MUV cell (0,1,1,2) to the all-2 cell, and (c) a green 16-cell PI loop extending from the purple (1, 1, 1, 1)-type MUV cell (1,1,1,1) to the all-2 cell. The total number of PI loops is c(4, 2) + 2 c(4, 2) + c(4, 4) = 6 + 12 + 1 = 19. The binary success at level 2 is given by the minimal sum-of-products expression:

$$S\{\geq 2\} = S\{2,3\} = X_1\{\geq 1\} X_2\{\geq 1\} X_3\{\geq 1\} X_4\{\geq 1\} \vee X_1\{\geq 0\} X_2\{\geq 1\} X_3\{\geq 1\} X_4\{\geq 2\} \vee X_1\{\geq 2\} X_2\{\geq 1\} X_3\{\geq 0\} X_4\{\geq 1\} \vee X_1\{\geq 1\} X_2\{\geq 0\} X_3\{\geq 2\} X_4\{\geq 1\} \vee X_1\{\geq 2\} X_2\{\geq 1\} X_3\{\geq 1\} X_4\{\geq 0\} \vee X_1\{\geq 0\} X_2\{\geq 1\} X_3\{\geq 2\} X_4\{\geq 1\} \vee X_1\{\geq 1\} X_2\{\geq$$

$1\} X_3\{\geq 2\} X_4\{\geq 0\} \vee X_1\{\geq 2\} X_2\{\geq 0\} X_3\{\geq 1\} X_4\{\geq 1\} \vee X_1\{\geq 1\} X_2\{\geq 2\} X_3\{\geq 0\} X_4\{\geq 1\} \vee X_1\{\geq 0\} X_2\{\geq 2\} X_3\{\geq 1\} X_4\{\geq 1\} \vee X_1\{\geq 1\} X_2\{\geq 0\} X_3\{\geq 1\} X_4\{\geq 2\} \vee X_1\{\geq 1\} X_2\{\geq 1\} X_3\{\geq 0\} X_4\{\geq 2\} \vee X_1\{\geq 1\} X_2\{\geq 2\} X_3\{\geq 1\} X_4\{\geq 0\} \vee X_1\{\geq 0\} X_2\{\geq 0\} X_3\{\geq 2\} X_4\{\geq 2\} \vee X_1\{\geq 0\} X_2\{\geq 2\} X_3\{\geq 0\} X_4\{\geq 2\} \vee X_1\{\geq 0\} X_2\{\geq 2\} X_3\{\geq 2\} X_4\{\geq 0\} \vee X_1\{2\} X_2\{0\} X_3\{0\} X_4\{2\} \vee X_1\{\geq 2\} X_2\{\geq 0\} X_3\{\geq 2\} X_4\{\geq 0\} \vee X_1\{\geq 2\} X_2\{\geq 2\} X_3\{\geq 0\} X_4\{\geq 0\}.$ (29)

Figure 7c covers the binary success at level 2 with 19 disjoint (non-overlapping) loops, comprising: (a) a single intact 16-cell loop, whose MUV is the sole cell of the (1, 1, 1, 1)-type, (b) 12 loops (four curtailed 4-cell loops, four curtailed 2-cell loops, and four curtailed 1-cell loops) that inherit the ones whose MUVs are (2, 1, 1, 0)-type, and (c) six loops (each a curtailed 1-cell loop containing just the MUV) that inherit the ones whose MUVs are (2, 2, 0, 0)-type. This threshold function is a shellable one, with each of the original PI loops retained as a single disjoint loop. The PRE for the binary success at level 2 is:

$S_{PRE}\{\geq 2\} = S_{PRE}\{2,3\} = X_1\{1,2\} X_2\{1,2\} X_3\{1,2\} X_4\{1,2\} \vee X_1\{0\} X_2\{1,2\} X_3\{1,2\} X_4\{2\} \vee X_1\{2\} X_2\{1,2\} X_3\{0\} X_4\{1,2\} \vee X_1\{1,2\} X_2\{0\} X_3\{2\} X_4\{1,2\} \vee X_1\{2\} X_2\{1,2\} X_3\{1,2\} X_4\{0\} \vee X_1\{0\} X_2\{1,2\} X_3\{2\} X_4\{1\} \vee X_1\{1\} X_2\{1,2\} X_3\{2\} X_4\{0\} \vee X_1\{2\} X_2\{0\} X_3\{1\} X_4\{1,2\} \vee X_1\{1\} X_2\{2\} X_3\{0\} X_4\{1,2\} \vee X_1\{0\} X_2\{2\} X_3\{1\} X_4\{1\} \vee X_1\{1\} X_2\{0\} X_3\{1\} X_4\{2\} \vee X_1\{1\} X_2\{1\} X_3\{0\} X_4\{2\} \vee X_1\{1\} X_2\{2\} X_3\{1\} X_4\{0\} \vee X_1\{0\} X_2\{0\} X_3\{2\} X_4\{2\} \vee X_1\{0\} X_2\{2\} X_3\{0\} X_4\{2\} \vee X_1\{0\} X_2\{2\} X_3\{2\} X_4\{0\} \vee X_1\{2\} X_2\{0\} X_3\{0\} X_4\{2\} \vee X_1\{2\} X_2\{0\} X_3\{2\} X_4\{0\} \vee X_1\{2\} X_2\{2\} X_3\{0\} X_4\{0\}.$ (30)

Variants of Fig. 8 deal with binary success at level 1, i.e., with $S\{\geq 1\} = S\{1, 2, 3\}$. Figure 8a shows the binary success at level 1 with MUV cells highlighted in yellow for the (2, 0, 0, 0)-type and in green for the (1, 1, 0, 0)-type. Figure 8b demonstrates two sample prime-implicant (PI) loops for the binary success at level 1. These include a 27-cell PI loop extending over 3 rows from the yellow (2, 0, 0, 0)-type MUV cell (0,0,0,2) to the all-2 cell (2, 2, 2, 2), and a 36-cell PI loop extending over 4 columns from the green (1, 1, 0, 0)-type MUV cell (1,1, 0, 0) to the all-2 cell. The total number of PI loops is c(4, 1) + c(4, 2) = 4 + 6 = 10. The binary success at level 1 is given by the minimal sum-of-products expression:

$S\{\geq 1\} = S\{1,2,3\} = X_1\{\geq 1\} X_2\{\geq 0\} X_3\{\geq 1\} X_4\{\geq 0\} \vee X_1\{\geq 0\} X_2\{\geq 1\} X_3\{\geq 1\} X_4\{\geq 0\} \vee X_1\{\geq 1\} X_2\{\geq 0\} X_3\{\geq 0\} X_4\{\geq 1\} \vee X_1\{\geq 0\} X_2\{\geq 1\} X_3\{\geq 0\} X_4\{\geq 1\} \vee X_1\{\geq 1\} X_2\{\geq 1\} X_3\{\geq 0\} X_4\{\geq 0\} \vee X_1\{\geq 0\} X_2\{\geq 0\} X_3\{\geq 1\} X_4\{\geq 1\} \vee X_1\{\geq 0\} X_2\{\geq 0\} X_3\{\geq 0\} X_4\{\geq 2\} \vee X_1\{\geq 0\} X_2\{\geq 0\} X_3\{\geq 2\} X_4\{\geq 0\} \vee X_1\{\geq 0\} X_2\{\geq 2\} X_3\{\geq 0\} X_4\{\geq 0\} \vee X_1\{\geq 2\} X_2\{\geq 0\} X_3\{\geq 0\} X_4\{\geq 0\}.$ (31)

Figure 8c covers the binary success at level 1 with ten disjoint (non-overlapping) loops, comprising: (a) six loops (an intact 36-cell loop and two curtailed 12-cell loops plus three curtailed 4-cell loops) that inherit the ones whose MUVs are (1, 1, 0, 0)-type, and (b) four loops (each a curtailed 1-cell loop containing just the yellow MUV) that inherit the ones whose MUVs are (2, 0, 0, 0)-type. This threshold function is a shellable one, with each of the original PI loops retained as a single disjoint loop. The PRE for the binary success at level 1 is:

$S_{PRE}\{\geq 1\} = S_{PRE}\{1,2,3\} = X_1\{1,2\} X_2\{0,1,2\} X_3\{1,2\} X_4\{0,1,2\} \vee X_1\{0\} X_2\{1,2\} X_3\{1,2\} X_4\{0,1,2\} \vee X_1\{1,2\} X_2\{0,1,2\} X_3\{0\} X_4\{1,2\} \vee X_1\{0\} X_2\{1,2\} X_3\{0\} X_4\{1,2\} \vee X_1\{1,2\} X_2\{1,2\} X_3\{0\} X_4\{0\} \vee$

$$X_1\{0\}\, X_2\{0\}\, X_3\{1,2\}\, X_4\{1,2\} \vee X_1\{0\}\, X_2\{0\}\, X_3\{0\}\, X_4\{2\} \vee X_1\{0\}\, X_2\{0\}\, X_3\{2\}\, X_4\{0\} \vee$$
$$X_1\{0\}\, X_2\{2\}\, X_3\{0\}\, X_4\{0\} \vee X_1\{2\}\, X_2\{0\}\, X_3\{0\}\, X_4\{0\}. \qquad (32)$$

We obtained equations (26)-(32) through appropriate reading of the MVKMs in variants of Figs. 6-8. Alternatively, we could have derived the minimal sum-of-products expressions (26). (29), and (31) through simple inspection of the MUVs in Table 4. We could also have obtained the PREs in equations (27). (28), (30) and (32) through the purely algebraic techniques of multi-state disjointing (see Appendix A) or multi-state Boole-Shannon Expansion (see Appendix B). The PREs obtained can be readily converted (on a one-to-one basis) into expected values by replacing the logical ORing and ANDing by arithmetic counterparts of addition and multiplication and replacing component instances by their expected values. For example, we can rewrite (28) as:

$$E\{S_{PRE}\{\geq 3\}\} = E\{S\{\geq 3\}\} = E\{X_1\{2\}\}\, E\{X_2\{1,2\}\}\, E\{X_3\{2\}\}\, E\{X_4\{1,2\}\} +$$
$$E\{X_1\{2\}\}\, E\{X_2\{1,2\}\}\, E\{X_3\{1\}\}\, E\{X_4\{2\}\} + E\{X_1\{1\}\}\, E\{X_2\{1,2\}\}\, E\{X_3\{2\}\}\, E\{X_4\{2\}\} +$$
$$E\{X_1\{2\}\}\, E\{X_2\{2\}\}\, E\{X_3\{1\}\}\, E\{X_4\{1\}\} + E\{X_1\{1\}\}\, E\{X_2\{2\}\}\, E\{X_3\{1\}\}\, E\{X_4\{2\}\} +$$
$$E\{X_1\{1\}\}\, E\{X_2\{2\}\}\, E\{X_3\{2\}\}\, E\{X_4\{1\}\} + E\{X_1\{2\}\}\, E\{X_2\{2\}\}\, E\{X_3\{2\}\}\, E\{X_4\{0\}\} +$$
$$E\{X_1\{2\}\}\, E\{X_2\{2\}\}\, E\{X_3\{0\}\}\, E\{X_4\{2\}\} + E\{X_1\{2\}\}\, E\{X_2\{0\}\}\, E\{X_3\{2\}\}\, E\{X_4\{2\}\} +$$
$$E\{X_1\{0\}\}\, E\{X_2\{2\}\}\, E\{X_3\{2\}\}\, \{X_4\{2\}\}. \qquad (33)$$

Once the upper-value (level success) expectations $E\{S\{\geq 3\}\}, E\{S\{\geq 2\}\}$, and $E\{S\{\geq 1\}\}$ are obtained, the instance expectations are computed as:

$$E\{S\{3\}\} = E\{S\{\geq 3\}\}, \qquad (34a)$$

$$E\{S\{2\}\} = E\{S\{\geq 2\}\} - E\{S\{\geq 3\}\}, \qquad (34b)$$

$$E\{S\{1\}\} = E\{S\{\geq 1\}\} - E\{S\{\geq 2\}\}, \qquad (34c)$$

$$E\{S\{0\}\} = 1.0 - E\{S\{\geq 1\}\}. \qquad (34d)$$

Having completed our solution in a success perspective, we now proceed to present a dual solution in a failure perspective. However, we omit the detailed discussion of minimal expressions, and go directly to the PRE ones.

Variants of Fig. 9 deal with binary failure at level 3, i.e., with $S\{\leq 2\} = S\{< 3\} = S\{0,1,2\}$. Figure 9a highlights the MLV cells for this binary failure in green for the (2, 2, 1, 0)-type, and in yellow for the (2, 1, 1, 1)-type. Figure 9b covers the binary failure at level 3 with 16 disjoint (non-overlapping) loops, comprising: (a) four loops (an intact 24-cell loop plus three curtailed 8-cell loops) that inherit the ones whose MLVs are of the (2, 1, 1, 1)-type, and (b) 12 loops (6 curtailed 2-cell loops plus 6 curtailed 1-cell loops (just the MLVs)) that inherit the ones whose MLVs are (2, 2, 1, 0)-type. Once more, this threshold function is a shellable one, with each of the original PI loops retained as a single disjoint loop. The PRE for the binary failure at level 3 is:

$$S_{PRE}\{\leq 2\} = S_{PRE}\{< 3\} = S_{PRE}\{0,1,2\} = X_1\{0,1\}\, X_2\{0,1,2\}\, X_3\{0,1\}\, X_4\{0,1\} \vee$$
$$X_1\{2\}\, X_2\{0,1\}\, X_3\{0,1\}\, X_4\{0,1\} \vee X_1\{0,1\}\, X_2\{0,1\}\, X_3\{2\}\, X_4\{0,1\} \vee$$
$$X_1\{0,1\}\, X_2\{0,1\}\, X_3\{0,1\}\, X_4\{2\} \vee X_1\{0\}\, X_2\{0,1\}\, X_3\{2\}\, X_4\{2\} \vee$$
$$X_1\{2\}\, X_2\{2\}\, X_3\{0\}\, X_4\{0,1\} \vee X_1\{0\}\, X_2\{2\}\, X_3\{2\}\, X_4\{0,1\} \vee X_1\{0\}\, X_2\{2\}\, X_3\{0,1\}\, X_4\{2\} \vee$$
$$X_1\{2\}\, X_2\{0,1\}\, X_3\{0\}\, X_4\{2\} \vee X_1\{2\}\, X_2\{0,1\}\, X_3\{2\}\, X_4\{0\} \vee X_1\{1\}\, X_2\{0\}\, X_3\{2\}\, X_4\{2\} \vee$$
$$X_1\{1\}\, X_2\{2\}\, X_3\{0\}\, X_4\{2\} \vee X_1\{1\}\, X_2\{2\}\, X_3\{2\}\, X_4\{0\} \vee X_1\{2\}\, X_2\{0\}\, X_3\{1\}\, X_4\{2\} \vee$$
$$X_1\{2\}\, X_2\{0\}\, X_3\{2\}\, X_4\{1\} \vee X_1\{2\}\, X_2\{2\}\, X_3\{1\}\, X_4\{0\}. \qquad (35)$$

Variants of Fig. 10 deal with binary failure at level 2, i.e., with $S\{\leq 1\} = S\{< 2\} = S\{0, 1\}$. Figure 10a highlights the MLV cells of this binary failure in green for the (2, 1, 0, 0)-type, and in yellow for the (1, 1, 1, 0)-type. Figure 10b covers the binary failure at level 2 with 16 disjoint (non-overlapping) loops, comprising: (a) four loops (an intact 8-cell loop plus three curtailed 4-cell, 2-cell and 1-cell loops) that inherit the ones whose MLVs are (1, 1, 1, 0)-type, and (b) 12 loops (4 curtailed 2-cell loops plus 8 curtailed 1-cell loops (just the MLVs)) that inherit the ones whose MLVs are (2, 1, 0, 0)-type. This threshold function is a shellable one, with each of the original PI loops retained as a single disjoint loop. The PRE for the binary failure at level 2 is:

$$S_{PRE}\{\leq 1\} = S_{PRE}\{< 2\} = S_{PRE}\{0, 1\} = X_1\{0\} X_2\{0, 1\} X_3\{0, 1\} X_4\{0, 1\} \vee$$
$$X_1\{1\} X_2\{0, 1\} X_3\{0\} X_4\{0, 1\} \vee X_1\{1\} X_2\{0, 1\} X_3\{1\} X_4\{0\} \vee X_1\{1\} X_2\{0\} X_3\{1\} X_4\{1\} \vee$$
$$X_1\{0\} X_2\{0, 1\} X_3\{0\} X_4\{2\} \vee X_1\{0\} X_2\{0, 1\} X_3\{2\} X_4\{0\} \vee X_1\{0\} X_2\{2\} X_3\{0\} X_4\{0, 1\} \vee$$
$$X_1\{2\} X_2\{0\} X_3\{0\} X_4\{0, 1\} \vee X_1\{0\} X_2\{0\} X_3\{1\} X_4\{2\} \vee X_1\{0\} X_2\{0\} X_3\{2\} X_4\{1\} \vee$$
$$X_1\{0\} X_2\{2\} X_3\{1\} X_4\{0\} \vee X_1\{1\} X_2\{0\} X_3\{0\} X_4\{2\} \vee X_1\{1\} X_2\{0\} X_3\{2\} X_4\{0\} \vee$$
$$X_1\{1\} X_2\{2\} X_3\{0\} X_4\{0\} \vee X_1\{2\} X_2\{0\} X_3\{1\} X_4\{0\} \vee X_1\{2\} X_2\{1\} X_3\{0\} X_4\{0\}. \quad (36)$$

Variants of Fig. 11 deal with binary failure at level 2, i.e., with $S\{\leq 0\} = S\{< 1\} = S\{0\}$. Figure 11a highlights the MLV cells of the binary failure at level 1 in yellow for its sole (1, 0, 0, 0)-type. Figure 11b covers this binary failure with four disjoint (non-overlapping) loops, comprising an intact 2-cell loop plus three curtailed 1-cell loops that inherit the ones whose MLVs are (1, 0, 0, 0)-type. The PRE for the binary failure at level 1 is:

$$S_{PRE}\{\leq 0\} = S_{PRE}\{< 1\} = S_{PRE}\{0\} = X_1\{0\} X_2\{0, 1\} X_3\{0\} X_4\{0\} \vee$$
$$X_1\{0\} X_2\{0\} X_3\{0\} X_4\{1\} \vee X_1\{0\} X_2\{0\} X_3\{1\} X_4\{0\} \vee X_1\{1\} X_2\{0\} X_3\{0\} X_4\{0\}. \quad (37)$$

Once the lower-value (level failure) expectations $E\{S\{\leq 0\}\}, E\{S\{\leq 1\}\}$, and $E\{S\{\leq 2\}\}$ are obtained, the instance expectations are computed as:

$$E\{S\{0\}\} = E\{S\{\leq 0\}\}, \quad (38a)$$

$$E\{S\{1\}\} = E\{S\{\leq 1\}\} - E\{S\{\leq 0\}\}, \quad (38b)$$

$$E\{S\{2\}\} = E\{S\{\leq 2\}\} - E\{S\{\leq 1\}\}, \quad (38c)$$

$$E\{S\{3\}\} = 1.0 - E\{S\{\leq 2\}\}. \quad (38d)$$

## 9. Conclusions

This paper utilizes map tools for the reliability characterization and analysis of general threshold multi-state coherent systems, which are non-repairable systems with independent non-identical components. The paper presents switching-algebraic expressions of both system success and system failure at each non-zero level. These expressions are given as minimal sum-of-products formulas for system success or as probability–ready expressions for both system success and failure. The paper utilizes a convenient map representation via the multi-valued Karnaugh map for the system structure function $S$, or via $M$ maps of binary entries and multi-valued inputs representing the success/failure at every non-zero level of the system. Further system characterizations are also given in terms of minimal upper vectors or maximal lower vectors. Our results demonstrate that, similarly to binary threshold systems, multi-state threshold systems possess the shellability property. In fact, the minimal sum-of-products formulas for system success or system failure at each non-zero level can be converted to disjoint formulas without any change in the number of products. Great emphasis is placed on making a minimal departure from binary concepts and techniques, while taking care to clarify novel issues that emerged due to generalizations introduced by the multi-state model.

The binary threshold model [23-32] has been developed initially as a model of system reliability. However, it might find further applications in other disciplines. One notable such discipline is the one of weighted voting systems [40-44]. These systems have been traditionally treated as binary ones, but their extension to multi-state ones has recently witnessed rapid progress. Such an extension is useful for the treatment of weighted voting systems in which a voter does not necessarily cast a yes-no vote but might abstain or be absent as well [77-79]. Another notable discipline is one in which the topics of system reliability and weighted voting systems are combined [80-82].

**Appendix A. Multi-State Probability-Ready Expressions**

The concept of a probability-ready expression (RRE) is of great utility in the two-valued logical domain [26, 48-53], and it is still useful for the multi-valued logical domain [45, 56, 59, 61-66]. A Probability-Ready Expression is a random expression that can be directly transformed, on a one-to-one basis, to its statistical expectation (its probability of being equal to 1) by the replacement of all logic variables by their statistical expectations, with the substitution of logical multiplication and addition (ANDing and ORing) by their arithmetic counterparts. A logic expression is a PRE if

   a)  all ORed products (terms formed by ANDing) are disjoint (mutually exclusive), and
   b)  all ANDed sums (alterms formed via ORing) are statistically independent.

Condition (a) is satisfied if for every pair of ORed terms, there is at least a single opposition, i.e., there is at least one variable that appears with a certain set of instances in one term and appears with a complementary set of instances in the other.  Condition (b) is satisfied if for every pair of ANDed alterms (sums of disjunctions of literals), one alterm involves variables describing a certain set of components, while the other alterm depends on variables describing a set of different components (under the assumption of independence of components).

While there are many methods to introduce characteristic (a) of orthogonality (disjointness) into a Boolean expression, there is no way to induce characteristic (b) of statistical independence. The best that one can do is to observe statistical independence when it exists, and then be careful not to destroy or spoil it but stive to take advantage of it. Since one has the freedom of handling a problem from a success or a failure perspective, a choice should be made as to which of the two perspectives can more readily produce a PRE form.

The introduction of orthogonality might be achieved as follows. If neither of the two terms $A$ and $B$ in the sum $(A \vee B)$ subsumes the other $(A \vee B \neq A$ and $A \vee B \neq B)$ and the two terms are not disjoint $(A \wedge B \neq 0)$, then $B$ can be disjointed with $A$ by the factorization of any common factor and then applying the *Reflection Law*, namely:

$$A \vee B = C((A/C) \vee (B/C)) = C((A/C) \vee \overline{(A/C)}(B/C)) = A \vee \overline{(A/C)}B. \quad (A.1)$$

In (A.1), the symbol $C$ denotes the common factor of $A\ and\ B$, and the multi-valued quotient $(A/C)$ might be viewed as the term $A$ with its part common with $B$ removed or set to 1. Note that $(A.1)$ leaves the term $A$ intact and replaces the term $B$ by an expression that is disjoint with A. The quotient $(A/C)$ is a product of $e$ entities  $Y_k$ $(1 \leq$

$k \leq e$), so that $\overline{(A/C)}$ might be expressed via *De Morgan's Law* as a disjunction of the form:

$$\overline{(A/C)} = \bigvee_{k=1}^{e} \bar{Y}_k. \qquad (A.2)$$

Note that each $Y_K$ stands for a disjunction of certain instances of some variable $X_{i(k)}$ and hence $\bar{Y}_k$ is a disjunction of the complementary instances of the same variable. If we combine (A.1) with (A.2), we realize that the term $B$ is replaced by $e$ terms ($e \geq 1$), which are each disjoint with the term $A$, but are not necessarily disjoint among themselves. Therefore, we replace the *De Morgan's Law* in (A.2) by a disjoint version of it [26], namely:

$$\overline{(A/C)} = \bar{Y}_1 \vee Y_1 \bar{Y}_2 \vee Y_1 Y_2 \bar{Y}_3 \vee \ldots \vee Y_1 Y_2 \ldots Y_{e-1} \bar{Y}_e$$
$$= \bar{Y}_1 \vee Y_1 (\bar{Y}_2 \vee Y_2 (\bar{Y}_3 \vee \ldots \ldots \vee (\bar{Y}_{e-1} \vee Y_{e-1} \bar{Y}_e) \ldots)). \qquad (A.3)$$

When (26) is combined with (24), the first term $A$ remains intact, while the second term $B$ is replaced by $e$ terms which are each disjoint with $A$ and are also disjoint among themselves. Note that one has a choice of either disjointing $B$ with $A$ in $A \vee B$, or disjointing $A$ with $B$ in $B \vee A$. A greedy practice that is likely to yield good results is to order the terms in a given disjunction so that those with fewer literals appear earlier.

### Appendix B. The Multi-Valued Boole-Shannon Expansion

The most effective way for converting a Boolean formula into a PRE form is the Boole-Shannon Expansion, which takes the following form in the two-valued case [26, 50, 68, 70, 76]

$$f(\mathbf{X}) = (\bar{X}_i \wedge f(\mathbf{X}|0_i)) \vee (X_i \wedge f(\mathbf{X}|1_i)), \qquad (B.1)$$

This Boole-Shannon Expansion expresses a (two-valued) Boolean function $f(\mathbf{X})$ in terms of its two subfunctions $f(\mathbf{X}|0_i)$ and $f(\mathbf{X}|1_i)$. These subfunctions are equal to the Boolean quotients $f(\mathbf{X})/\bar{X}_i$ and $f(\mathbf{X})/X_i$, and hence are obtained by restricting $X_i$ in the expression of $f(\mathbf{X})$ to 0 and 1, respectively. If $f(\mathbf{X})$ is a sum-of-products (sop) expression of $n$ variables, the two sub-functions $f(\mathbf{X}|0_i)$ and $f(\mathbf{X}|1_i)$ are functions of at most $(n-1)$ variables. A multi-valued extension of (B.1) is

$$S(\mathbf{X}) = X_i\{0\} \wedge (S(\mathbf{X})/X_i\{0\}) \vee X_i\{1\} \wedge (S(\mathbf{X})/X_i\{1\}) \vee X_i\{2\} \wedge (S(\mathbf{X})/X_i\{2\}) \vee X_i\{3\} \wedge (S(\mathbf{X})/X_i\{3\}) \vee \ldots \vee X_i\{m_i\} \wedge (S(\mathbf{X})/X_i\{m_i\}). \qquad (B.2)$$

Rushdi and Ghaleb [59] provide a formal proof of (B.2), and note that once the sub-functions in (B.2) are expressed by PRE expressions, $S(\mathbf{X})$ will be also in PRE form. The expansion (B.2) might be viewed as a justification of the construction of the multi-valued Karnaugh map used extensively herein.

**Conflict of Interest**

The authors assert that no conflict of interest exists.

**Acknowledgment**

The first named author (AMAR) is greatly indebted to Engineer Mahmoud Ali Rushdi of Munich, Germany for fruitful discussions and continuous collaborations.

# References

[ 1] **Barlow, R. E.,** and **Wu, A. S**., Coherent systems with multi-state components, *Mathematics of Operations Research*, **3**(4): 275-281 (1978).

[ 2] **El-Neweihi, E., Proschan, F.,** and **Sethuraman, J**., Multistate coherent systems, *Journal of Applied Probability*, **15**(4): 675-688 (1978).

[ 3] **Griffith, W. S.,** Multistate reliability models, *Journal of Applied Probability*, **17(**3): 735-744 (1980).

[ 4] **Fardis, M. N.,** and **Cornell, C. A**., Analysis of coherent multistate systems, *IEEE Transactions on Reliability*, **30**(2): 117-122 (1981).

[ 5] **Fardis, M. N.,** and **Conrnell, C. A**., Multistate reliability analysis, *Nuclear Engineering and Design*, **71**(3): 329-336 (1982).

[ 6] **Hudson, J. C.,** and **Kapur, K. C.,** Reliability analysis for multistate systems with multistate components, *AIIE Transactions*, **15**(2): 127-135 (1983).

[ 7] **El-Neweihi, E.,** and **Proschan, F.**, Degradable systems: a survey of multistate system theory, *Communications in Statistics-Theory and Methods*, **13**(4): 405-432 (1984).

[ 8] **Wood, A. P.,** Multistate block diagrams and fault trees. *IEE, Transactions on Reliability*, **34**(3): 236-240 (1985).

[ 9] **Janan, X.,** On multistate system analysis, *IEEE Transactions on Reliability*, **34**(4): 329-337 (1985).

[ 10] **Hudson, J. C.,** and **Kapur, K. C.** Reliability bounds for multistate systems with multistate components, *Operations Research*, **33**(1), 153-160 (1985).

[ 11] **Boedigheimer, R. A.,** and **Kapur, K. C.,** Customer-driven reliability models for multistate coherent systems. *IEEE Transactions on Reliability*, **43**(1): 46-50 (1994).

[ 12] **Xue, J.,** and **Yang, K.,** Symmetric relations in multistate systems. *IEEE Transactions on Reliability*, **44**(4): 689-693 (1995).

[ 13] **Levitin, G., Lisnianski, A**., and **Ushakov, I.,** Reliability of multi-state systems: a historical overview. In **Lindqvist, B. H.** and **Doksum, K. A.** (Editors), *Mathematical and Statistical Methods in Reliability*, World Scientific: pp. 123-137 (2003).

[ 14] **Huang, J.,** and **Zuo, M. J.,** Dominant multi-state systems. *IEEE Transactions on Reliability*, **53**(3): 362-368 (2004).

[ 15] **Zuo, M. J., Tian, Z.,** and **Huang, H. Z**., An efficient method for reliability evaluation of multistate networks given all minimal path vectors, *IIE Transactions*, **39**(8): 811-817 (2007).

[ 16] **Xing, L.,** and **Dai, Y. S**., A new decision-diagram-based method for efficient analysis on multistate systems, *IEEE Transactions on Dependable and Secure Computing*, **6**(3): 161-174 (2008).

[ 17] **Yeh, W. C., Lin L. E.** and **Chen Y. C**., Reliability evaluation of multi-state quick path

flow networks, *Journal of Quality*, **1**;20:199–215 (2013).

[ 18] **Lin, Y. K., Huang, C. F**., and **Yeh, C. T.,** Network reliability with deteriorating product and production capacity through a multi-state delivery network. International Journal of Production Research, **52**(22): 6681-6694 (2014).

[ 19] **Yeh, W**. **C.**, An improved sum-of-disjoint-products technique for symbolic multi-state flow network reliability, *IEEE Transactions on Reliability*, **64**(4): 1185–1193 (2015).

[ 20] **Bai, G., Tian, Z.,** and **Zuo, M. J**., Reliability evaluation of multistate networks: An improved algorithm using state-space decomposition and experimental comparison. *IISE Transactions*, **50**(5): 407–418 (2018).

[ 21] **Lin Y. K., Nguyen T. P.** and **Yeng L. C. L**., Reliability evaluation of a multi-state air transportation network meeting multiple travel demands. *Annals of Operations Research*, **277**(1): 63–82 (2019).

[ 22] **Negi, M., Shah, M., Kumar, A., Ram, M.,** and **Saini, S.** Assessment of Reliability Function and Signature of Energy Plant Complex System. In *Reliability and Maintainability Assessment of Industrial Systems* (pp. 257-269). Springer, Cham (2022).

[ 23] **Rushdi, A. M.** Threshold systems and their reliability, *Microelectronics and Reliability*, **30**(2), 299-312 (1990).

[ 24] **Rushdi, A. M. A**., and **Alturki, A. M**., Reliability of coherent threshold systems, *Journal of Applied Sciences*, **15**(3), 431-443 (2015).

[ 25] **Rushdi, A. M. A.,** and **Bjaili, H. A.**, An ROBDD algorithm for the reliability of double-threshold systems. *British Journal of Mathematics and Computer Science*, **19**(6), 1-17 (2016).

[ 26] **Rushdi, A. M.,** and **Rushdi, M. A.**, Switching-Algebraic Analysis of System Reliability. Chapter 6 in **Ram, M.** and **Davim, P**. (Editors), *Advances in Reliability and System Engineering*, Management and Industrial Engineering Series, Springer International Publishing, Cham, Switzerland: 139-161 (2017).

[ 27] **Rushdi, A. M. A**., and **Alturki, A. M**., Novel representations for a coherent threshold reliability system: a tale of eight signal ow graphs, *Turkish Journal of Electrical Engineering & Computer Sciences*, **26(**1), 257-269 (2018).

[ 28] **Rushdi, A. M.** and **Alturki, A. M**., Representations of a coherent reliability system via signal flow graphs, *Journal of King Abdulaziz University: Engineering Sciences*, **31**(1): 3-17 (2020).

[ 29] **Uswarman, R.,** and **Rushdi, A. M.** Reliability evaluation of rooftop solar photovoltaic using coherent threshold systems, *Journal of Engineering Research and Reports*, **20**(2), 32-44 (2021).

[ 30] **Muktiadji, R. F.,** and **Rushdi, A. M.** Reliability analysis of boost converters connected to a solar panel using a Markov approach. *Journal of Energy Research and Reviews*, **7**(1), 29-42 (2021).

[ 31] **Hidayat, T.,** and **Rushdi, A. M**. Reliability analysis of a home-scale microgrid based on a threshold system, *Journal of Energy Research and Reviews*, **7**(3), 14-26

(2021).

[32] **Budiman, F. N.** and **Rushdi, A. M**., Reliability evaluation of an electric power system utilizing fault-tree analysis, *Journal of Qassim University: Engineering and Computer Sciences*, 14(1) ( 2021).

[33] **Eryilmaz, S.** Capacity loss and residual capacity in weighted k-out-of-n: G systems, *Reliability Engineering & System Safety*, **136**, 140-144 (2015).

[34] **Li, X., You, Y.,** and **Fang, R.** On weighted k-out-of-n systems with statistically dependent component lifetimes, *Probability in the Engineering and Informational Sciences*, **30**(4), 533-546 (2016).

[35] **Zhang, Y.** Optimal allocation of active redundancies in weighted k-out-of-n systems, *Statistics & Probability Letters*, **135**, 110-117 (2018).

[36] **Salehi, M., Shishebor, Z.,** and **Asadi, M**. On the reliability modeling of weighted k-out-of-n systems with randomly chosen components. *Metrika*, **82**(5), 589-605 (2019).

[37] **Devrim, Y.,** and **Eryilmaz, S.** Reliability-based evaluation of hybrid wind–solar energy system. *Proceedings of the Institution of Mechanical Engineers, Part O: Journal of Risk and Reliability*, **235**(1), 136-143 (2021).

[38] **Eryilmaz, S.,** and **Ucum, K. A.** The lost capacity by the weighted k-out-of-n system upon system failure. *Reliability Engineering & System Safety*, **216**, Article number 107914 (2021).

[39] **Triantafyllou, I. S.** Reliability structures consisting of weighted components: Synopsis and new advances. *Journal of Reliability and Statistical Studies*, **16**(1) 25-56 (2023).

[40] **Alturki, A. M**., and **Rushdi, A. M. A.** Weighted voting systems: A threshold-Boolean perspective. *Journal of Engineering Research*, **4**(1), 125-143 (2016).

[41] **Rushdi, A. M.,** and **Rushdi, M. A**. Switching-algebraic calculation of Banzhaf voting indices, *arXiv Preprint*, 1-18 (2023). https://doi.org/10.48550/arXiv.2302.09367. Also in *Journal of Computational and Cognitive Engineering* (JCCE), **2**(3) (2023).

[42] **Rushdi, A. M.** and **Rushdi, M. A.** Towards a switching-algebraic theory of weighted monotone voting systems: The case of Banzhaf voting indices, *arXiv Preprint*, 1-34 (2023), https://doi.org/10.48550/arXiv.2302.09367.

[43] **Rushdi, A. M.** and **Rushdi, M. A.** Towards a switching-algebraic theory of weighted voting systems: Exploring restrictions on coalition formation, *arXiv Preprint*, 1-45 (2023), https://doi.org/10.48550/arXiv.2306.13684.

[44] **Rushdi, A. M.,** and **Rushdi, M. A.** Boolean-based probability: the case of a vector-weighted voting system. OFSPREPRINTS, https://osf.io/emsrj/, 1-21 (2023). Also in *King Abdulaziz University Journal: Engineering Sciences,* **34**(2) (2023).

[45] **Rushdi, A. M. A.,** and **Ghaleb, F. A. M.,** Reliability characterization of binary-imaged multi-state coherent threshold systems, *International Journal of Mathematical, Engineering and Management Sciences (IJMEMS)*, **6**(1): 309-321 (2021).

[46] **Ansell, J. I.,** and **Bendell, A.** On alternative definitions of multistate coherent systems, *Optimization*, **18**(1): 119-136 (1987).

[47] **Bennetts, R. G.,** On the analysis of fault trees, *IEEE Transactions on Reliability*, **R-24** (3): 175-185, (1975).

[48] **Bennetts, R. G.,** Analysis of reliability block diagrams by Boolean techniques, *IEEE Transactions on Reliability,* **R-31**(2): 159-166, (1982).

[49] **Rushdi, A. M.,** and **Al-Khateeb, D. L.,** A review of methods for system reliability analysis: A Karnaugh-map perspective, *Proceedings of the First Saudi Engineering Conference, Jeddah, Saudi Arabia*, vol. **1**: pp. 57-95, (1983).

[50] **Rushdi, A. M.,** and **Goda, A. S.**, Symbolic reliability analysis via Shannon's expansion and statistical independence, *Microelectronics and Reliability*, **25**(6): 1041-1053 (1985).

[51] **Rushdi, A. M.,** and **AbdulGhani, A. A.,** A comparison between reliability analyses based primarily on disjointness or statistical independence: The case of the generalized INDRA network, *Microelectronics and Reliability*, **33**(7): 965-978 (1993).

[52] **Rushdi, A. M. A.,** and **Hassan A. K.,** Reliability of migration between habitat patches with heterogeneous ecological corridors*, Ecological Modelling*, **304**: 1-10, (2015).

[53] **Rushdi, A. M. A.,** and **Hassan, A. K.,** An exposition of system reliability analysis with an ecological perspective, *Ecological Indicators*, **63**: 282-295, (2016).

[54] **Rushdi, A. M. A., Hassan, A. K.,** and **Moinuddin, M**., System reliability analysis of small-cell deployment in heterogeneous cellular networks, *Telecommunication Systems*, **73**: 371-381 (2019).

[55] **Rushdi, A. M. A.,** and **Hassan, A. K.**, On the Interplay between Reliability and Ecology, Chapter 35 in **Misra, K. B**. (Editor), *Advanced Handbook of Performability Engineering.* Springer Science & Business Media, (2020).

[56] **Rushdi, A. M. A.,** Utilization of symmetric switching functions in the symbolic reliability analysis of multi-state k-out-of-n systems, *International Journal of Mathematical, Engineering and Management Sciences (IJMEMS)*, **4**(2): 306-326 (2019).

[57] **Rushdi, A. M A**., and **Alsayegh, A. B**., Reliability analysis of a commodity-supply multi-state system using the map method. *Journal of Advances in Mathematics and Computer Science*, **31**(2): 1-17 (2019).

[58] **Rushdi, A. M., AlHuthali, S. A., AlZahrani N. A.,** and **Alsayegh, A. B.,** Reliability Analysis of Binary-Imaged Generalized Multi-State k-out-of-n Systems, *International Journal of Computer Science and Network Security (IJCSNS),* **20**(9): 251-264 (2020).

[59] **Rushdi A. M.** and **Ghaleb FA.,** Boolean-based symbolic analysis for the reliability of coherent multi-state systems of heterogeneous components. *Journal of King Abdulaziz University: Computing and Information Technology* Sciences, **9**(2): 1-25 (2020).

[60] **Rushdi, A. M. A.,** and **Amashah, M. H.**, Conventional and improved inclusion-

exclusion derivations of symbolic expressions for the reliability of a multi-state network, *Asian Journal of Research in Computer Science*. **8**(1): 21-45 (2021).

[ 61] **Rushdi, A. M. A.,** and **Amashah, M. H**., Symbolic derivation of a probability-ready expression for the reliability analysis of a multi-state delivery network, *Journal of Advances in Mathematics and Computer Science*. **36**(2): 37–56 (2021).

[ 62] **Rushdi A. M. A.,** and **Amashah M. H.,** Symbolic reliability analysis of a multi-state network. In *IEEE Fourth National Computing Colleges Conference (4th NCCC)*, Taif, Kingdom of Saudi Arabia, 1-4 (2021), DOI: 10.1109/NCCC49330.2021.9428843 (2021).

[ 63] **Rushdi, A. M.,** and **Amashah, M. H.**, A liaison among inclusion-exclusion, probability-ready expressions and Boole-Shannon expansion for multi-state reliability, *Journal of King Abdulaziz University: Computing and Information Technology Sciences*, **10**(2): 1-17 (2021).

[ 64] **Rushdi, A. M.,** and **Amashah, M. H.**, Derivation of minimal cutsets from minimal pathsets for a multi-state system and utilization of both sets in checking reliability expressions, *Journal of Engineering Research and Reports*, **20**(8): 22-33 (2021).

[ 65] **Rushdi, A. M., and Amashah, M. H.**, Reliability Analysis of a Multi-State Delivery Network through the Symbolic Derivation of a Probability-Ready Expression. Chapter 12 in *Recent Advances in Mathematical Research and Computer Science Vol. 4*, 113-132 (2021).

[ 66] **Rushdi, A. M.,** and **Amashah, M. H.**, Complementation of Multi-State System Success to Obtain System Failure and Utilization of Both Boolean Functions in Checking Reliability Expressions. Chapter 11 in *Novel Perspectives of Engineering Research Vol. 3*, 141-154 (2021).

[ 67] **Rushdi, A. M. A**., and **Alsayegh, A. B**., Karnaugh-Map Analysis of a Commodity-Supply Multi-State Reliability System. Chapter 3 in *Theory and Practice of Mathematics and Computer Science Vol. 8,* 20-37 (2021).

[ 68] **Rushdi, A. M.,** and **Ghaleb, F. A.** Switching-Algebraic Symbolic Analysis of the Reliability of Non-Repairable Coherent Multistate Systems. Chapter 8 in **M. Ram** and **Liudong Xing** (Editors), *Mathematics for Reliability Engineering: Modern Concepts and Applications*, De Gruyter Series on the Applications of Mathematics in Engineering and Information Sciences, De Gruyter, Berlin, Germany, 131–150 (2022).

[ 69] **Perkowski, M., Marek-Sadowska, M., Jozwiak, L., Luba, T., Grygiel, S., Nowicka, M., Malvi, R., Wang, Z.** and **Zhang, J.S.** Decomposition of multiple-valued relations. In *Proceedings 1997 27th IEEE International Symposium on Multiple-Valued Logic*: 13-18 (1997).

[ 70] **Rushdi, A. M. A.,** and **Ghaleb, F. A. M.,** A tutorial exposition of semi-tensor products of matrices with a stress on their representation of Boolean functions, *Journal of King Abdulaziz University: Computing and Information Technology Sciences*, **5**(1): 3-30 (2016).

[ 71] **Rushdi, A. M. A.,** Utilization of Karnaugh maps in multi- value qualitative comparative analysis, *International Journal of Mathematical, Engineering and Management Sciences* (*IJMEMS*), **3**(1): 28-46 (2018).

[ 72] **Rushdi, R. A.,** and **Rushdi, A. M**., Karnaugh-map utility in medical studies: The case of Fetal Malnutrition. *International Journal of Mathematical, Engineering and Management Sciences (IJMEMS)*, **3**(3): 220-244 (2018).

[ 73] **Rushdi, A. M. A.,** and **Al-Amoudi, M. A.,** Reliability analysis of a multi-state system using multi-valued logic, *IOSR Journal of Electronics and Communication Engineering (IOSR-JECE)*, **14**(1): 1-10 (2019).

[ 74] **Rushdi, A. M. A.,** and **Alsalami, O. M.** Multistate reliability evaluation of communication networks via Multi-Valued Karnaugh Maps and exhaustive search. *Advanced Aspects of Engineering Research Vol. 1*, 114-136 (2021).

[ 75] **Lu, Z., Liu, G.,** and **Liao, R. A.** pseudo Karnaugh mapping approach for datasets imbalance. In E3S Web of Conferences, **236** (04006) *EDP Sciences*,1-5 (2021).

[ 76] **Brown, F. M.,** *Boolean Reasoning: The Logic of Boolean Equations*, Kluwer Academic Publishers, Boston, USA (1990).

[ 77] **Freixas, J.,** and **Zwicker, W. S.** Weighted voting, abstention, and multiple levels of approval. *Social choice and welfare*, **21**, 399-431 (2003).

[ 78] **Freixas, J.,** and **Zwicker, W.** Anonymous yes–no voting with abstention and multiple levels of approval. *Games and Economic Behavior*, **67**(2), 428-444 (2009).

[ 79] **Bolle, F.** Voting with abstention. *Journal of Public Economic Theory*, **24**(1), 30-57 (2022).

[ 80] **Xie, M.,** and **Pham, H.** Modeling the reliability of threshold weighted voting systems. *Reliability Engineering & System Safety*, **87**(1), 53-63 (2005).

[ 81] **Ramirez-Marquez, J. E.** Holistic reliability analysis of weighted voting systems from a multi-state perspective. *IIE Transactions*, **40**(2): 122-132 (2007).

[ 82] **Liu, Q.,** and **Zhang, H**. Weighted voting system with unreliable links. *IEEE Transactions on Reliability*, **66**(2), 339-350 (2017).

Table 1. Meaning of component (engine) states

| Component state $X_i$ | Meaning |
|---|---|
| 0 | Failed |
| 1 | Half-power |
| 2 | Full-power |

Table 2. Meaning of system (aircraft) states

| System state $j$ | Meaning |
|---|---|
| 0 | Crash landing |
| 1 | Emergency landing on a foamed runway |
| 2 | Emergency landing on a normal runway |
| 3 | Normal landing |

Table 3. Relation between the range of weighted sum of component successes to the value or level of system success.

| Level $j$ | Corresponding Range $[T_j, T_{j+1})$ | Actual possible integer values of the weighted-Sum | Value of S |
|---|---|---|---|
| 0 | [0, 2) | {0, 1} | $S\{0\}$ |
| 1 | [2, 4) | {2, 3} | $S\{1\}$ |
| 2 | [4, 6) | {4, 5} | $S\{2\}$ |
| 3 | [6, 9) | {6, 7, 8} | $S\{3\}$ |

| $X_1$ | 0 | | | 1 | | | 2 | | | | |
|---|---|---|---|---|---|---|---|---|---|---|---|
| $X_2$ | 0 | 1 | 2 | 0 | 1 | 2 | 0 | 1 | 2 | | |
| | 0 | 1 | 2 | 1 | 2 | 3 | 2 | 3 | 4 | 0 | |
| | 1 | 2 | 3 | 2 | 3 | 4 | 3 | 4 | 5 | 1 | 0 |
| | 2 | 3 | 4 | 3 | 4 | 5 | 4 | 5 | 6 | 2 | |
| | 1 | 2 | 3 | 2 | 3 | 4 | 3 | 4 | 5 | 0 | |
| | 2 | 3 | 4 | 3 | 4 | 5 | 4 | 5 | 6 | 1 | 1 |
| | 3 | 4 | 5 | 4 | 5 | 6 | 5 | 6 | 7 | 2 | |
| | 2 | 3 | 4 | 3 | 4 | 5 | 4 | 5 | 6 | 0 | |
| | 3 | 4 | 5 | 4 | 5 | 6 | 5 | 6 | 7 | 1 | 2 |
| | 4 | 5 | 6 | 5 | 6 | 7 | 6 | 7 | 8 | 2 | |
| | | | | | | | | | | $X_4$ | $X_3$ |

$$X_1 + X_2 + X_3 + X_4$$

Fig. 1. A MVKM representing the arithmetic sum of four three-valued variables, $X_1, X_2 = X_3$ and $X_4$.

| $X_1$ | 0 | | | 1 | | | 2 | | | | |
|---|---|---|---|---|---|---|---|---|---|---|---|
| $X_2$ | 0 | 1 | 2 | 0 | 1 | 2 | 0 | 1 | 2 | | |
| | 0 | 0 | 1 | 0 | 1 | 1 | 1 | 1 | 2 | 0 | |
| | 0 | 1 | 1 | 1 | 1 | 2 | 1 | 2 | 2 | 1 | 0 |
| | 1 | 1 | 2 | 1 | 2 | 2 | 2 | 2 | 3 | 2 | |
| | 0 | 1 | 1 | 1 | 1 | 2 | 1 | 2 | 2 | 0 | |
| | 1 | 1 | 2 | 1 | 2 | 2 | 2 | 2 | 3 | 1 | 1 |
| | 1 | 2 | 2 | 2 | 2 | 3 | 2 | 3 | 3 | 2 | |
| | 1 | 1 | 2 | 1 | 2 | 2 | 2 | 2 | 3 | 0 | |
| | 1 | 2 | 2 | 2 | 2 | 3 | 2 | 3 | 3 | 1 | 2 |
| | 2 | 2 | 3 | 2 | 3 | 3 | 3 | 3 | 3 | 2 | |
| | | | | | | | | | | $X_4$ | $X_3$ |

$$S(X)$$

Fig. 2. A MVKM representing the structure (success) function of the running example.

| $X_1$ | | 0 | | | 1 | | | 2 | | | | |
|---|---|---|---|---|---|---|---|---|---|---|---|---|
| $X_2$ | 0 | 1 | 2 | 0 | 1 | 2 | 0 | 1 | 2 | | | |
| | 0 | 0 | 1 | 0 | 1 | 1 | 1 | 1 | 2 | 0 | | |
| | 0 | 1 | 1 | 1 | 1 | 2 | 1 | 2 | 2 | 1 | 0 | |
| | 1 | 1 | 2 | 1 | 2 | 2 | 2 | 2 | 3 | 2 | | |
| | 0 | 1 | 1 | 1 | 1 | 2 | 1 | 2 | 2 | 0 | | |
| | 1 | 1 | 2 | 1 | 2 | 2 | 2 | 2 | 3 | 1 | 1 | |
| | 1 | 2 | 2 | 2 | 2 | 3 | 2 | 3 | 3 | 2 | | |
| | 1 | 1 | 2 | 1 | 2 | 2 | 2 | 2 | 3 | 0 | | |
| | 1 | 2 | 2 | 2 | 2 | 3 | 2 | 3 | 3 | 1 | 2 | |
| | 2 | 2 | 3 | 2 | **3** | 3 | 3 | 3 | 3 | 2 | | |
| | | | | | | | | | | $X_4$ | $X_3$ | |

$S(X)$

Fig. 3. Demonstration that the cell/vector $X = [1, 2, 2, 2]$ (highlighted in yellow) is an MUV for level 3 since $S(X) = 3 \geq 3$, while all the cells/vectors $Y < X$ (highlighted in green) are such that $S(Y) < 3$.

| $X_1$ | | 0 | | | 1 | | | 2 | | | | |
|---|---|---|---|---|---|---|---|---|---|---|---|---|
| $X_2$ | 0 | 1 | 2 | 0 | 1 | 2 | 0 | 1 | 2 | | | |
| | 0 | 0 | 1 | 0 | 1 | 1 | 1 | 1 | 2 | 0 | | |
| | 0 | 1 | 1 | 1 | 1 | 2 | 1 | 2 | 2 | 1 | 0 | |
| | 1 | 1 | 2 | 1 | 2 | 2 | 2 | 2 | 3 | 2 | | |
| | 0 | 1 | 1 | 1 | 1 | 2 | 1 | 2 | 2 | 0 | | |
| | 1 | 1 | 2 | 1 | 2 | 2 | 2 | **2** | 3 | 1 | 1 | |
| | 1 | 2 | 2 | 2 | 2 | 3 | 2 | 3 | 3 | 2 | | |
| | 1 | 1 | 2 | 1 | 2 | 2 | 2 | 2 | 3 | 0 | | |
| | 1 | 2 | 2 | 2 | 2 | 3 | 2 | 3 | 3 | 1 | 2 | |
| | 2 | 2 | 3 | 2 | 3 | 3 | 3 | 3 | 3 | 2 | | |
| | | | | | | | | | | $X_4$ | $X_3$ | |

$S(X)$

Fig. 4. Demonstration that the cell/vector $X = [2, 1, 1, 1]$ (highlighted in yellow) is an MLV for level 2 since $S(X) = 2 \leq 2$, while all the cells/vectors $Y > X$ (highlighted in green) are such that $S(Y) > 2$.

| $X_1$ |   | 0 |   |   | 1 |   |   | 2 |   |       |       |
|-------|---|---|---|---|---|---|---|---|---|-------|-------|
| $X_2$ | 0 | 1 | 2 | 0 | 1 | 2 | 0 | 1 | 2 |       |       |
|       | 0 | 0 | 1 | 0 | 1 | 1 | 1 | 1 | 2 | 0     |       |
|       | 0 | 1 | 1 | 1 | 1 | 2 | 1 | 2 | 2 | 1     | 0     |
|       | 1 | 1 | 2 | 1 | 2 | 2 | 2 | 2 | 3 | 2     |       |
|       | 0 | 1 | 1 | 1 | 1 | 2 | 1 | 2 | 2 | 0     |       |
|       | 1 | 1 | 2 | 1 | 2 | 2 | 2 | 2 | 3 | 1     | 1     |
|       | 1 | 2 | 2 | 2 | 2 | 3 | 2 | 3 | 3 | 2     |       |
|       | 1 | 1 | 2 | 1 | 2 | 2 | 2 | 2 | 3 | 0     |       |
|       | 1 | 2 | 2 | 2 | 2 | 3 | 2 | 3 | 3 | 1     | 2     |
|       | 2 | 2 | 3 | 2 | 3 | 3 | 3 | 3 | 3 | 2     |       |
|       |   |   |   |   |   |   |   |   |   | $X_4$ | $X_3$ |

$$S(X)$$

Fig. 5. A MVKM representing the structure (success) function of the running example. Monotonicity with respect to component 1 can be observed by dividing the map into three submaps each comprising one third of the total map and consisting of three map columns. These three submaps are separated by two vertical $X_1$ borders transitioning from $X_1 = 0$ to $X_1 = 1$, and from $X_1 = 1$ to $X_1 = 2$, respectively. They represent the regions of $X_1 = 0$, $X_1 = 1$, and $X_1 = 2$ (corresponding to $X_1\{0\} = 1, X_1\{1\} = 1$, and $X_1\{2\} = 1$), and they stand for the multivalued quotients $S(X)/X_1\{0\}$, $S(X)/X_1\{1\}$, and $S(X)/X_1\{2\}$, which are each functions of $X_2$, $X_3$, and $X_4$. Monotonicity w.r.t. component 1 is demonstrated by the relations

$$S(X)/X_1\{2\} \geq S(X)/X_1\{1\} \geq S(X)/X_1\{0\}.$$

Table 4. Listing of sample elements of the set of Minimal Upper Vectors (MUVs) $\boldsymbol{\theta}(k)$, $(k = 3, 2, 1)$ through solution of the equation $X_1 + X_2 + X_3 + X_4 = 2k$ for integers $X_1, X_2, X_3, X_4$ lying in $\{0, 1, 2\}$.

| $k$ | $2k$ | Sample Solution | Corresponding prime-implicant loop | No of cells in the loop | No of variants of the sample solution |
|---|---|---|---|---|---|
| 3 | 6 | $[2, 2, 2, 0]$ | $X_1\{\geq 2\}X_2\{\geq 2\}X_3\{\geq 2\}X_4\{\geq 0\}$ | $1 * 1 * 1 * 3 = 3$ | $\binom{4}{1} = 4$ |
|   |   | $[2, 2, 1, 1]$ | $X_1\{\geq 2\}X_2\{\geq 2\}X_3\{\geq 1\}X_4\{\geq 1\}$ | $1 * 1 * 2 * 2 = 4$ | $\binom{4}{2} = 6$ |
| 2 | 4 | $[2, 2, 0, 0]$ | $X_1\{\geq 2\}X_2\{\geq 2\}X_3\{\geq 0\}X_4\{\geq 0\}$ | $1 * 1 * 3 * 3 = 9$ | $\binom{4}{2} = 6$ |
|   |   | $[2, 1, 1, 0]$ | $X_1\{\geq 2\}X_2\{\geq 1\}X_3\{\geq 1\}X_4\{\geq 0\}$ | $1 * 2 * 2 * 3 = 12$ | $\binom{4}{2}\binom{2}{1} = 12$ |
|   |   | $[1, 1, 1, 1]$ | $X_1\{\geq 1\}X_2\{\geq 1\}X_3\{\geq 1\}X_4\{\geq 1\}$ | $2 * 2 * 2 * 2 = 16$ | $\binom{4}{4} = 1$ |
| 1 | 2 | $[2, 0, 0, 0]$ | $X_1\{\geq 2\}X_2\{\geq 0\}X_3\{\geq 0\}X_4\{\geq 0\}$ | $1 * 3 * 3 * 3 = 27$ | $\binom{4}{1} = 4$ |
|   |   | $[1, 1, 0, 0]$ | $X_1\{\geq 1\}X_2\{\geq 1\}X_3\{\geq 0\}X_4\{\geq 0\}$ | $2 * 2 * 3 * 3 = 36$ | $\binom{4}{2} = 6$ |

Table 5. Listing of sample elements of the set of Maximal Lower Vectors (MLVs) $\sigma(k)$, $(k = 2, 1, 0)$ through solution of the equation $X_1 + X_2 + X_3 + X_4 = 2k + 1$ for integers $X_1, X_2, X_3, X_4$ lying in $\{0, 1, 2\}$.

| $k$ | $2k+1$ | Sample Solution | Corresponding prime-implicant loop | No of cells in the loop | No of variants of the sample solution |
|---|---|---|---|---|---|
| 2 | 5 | $[2, 2, 1, 0]$ | $X_1\{\leq 2\}X_2\{\leq 2\}X_3\{\leq 1\}X_4\{\leq 0\}$ | $3*3*2*1 = 18$ | $\binom{4}{2}\binom{2}{1} = 12$ |
|   |   | $[2, 1, 1, 1]$ | $X_1\{\leq 2\}X_2\{\leq 1\}X_3\{\leq 1\}X_4\{\leq 1\}$ | $3*2*2*2 = 24$ | $\binom{4}{1} = 4$ |
| 1 | 3 | $[2, 1, 0, 0]$ | $X_1\{\leq 2\}X_2\{\leq 1\}X_3\{\leq 0\}X_4\{\leq 0\}$ | $3*2*1*1 = 6$ | $\binom{4}{2}\binom{2}{1} = 12$ |
|   |   | $[1, 1, 1, 0]$ | $X_1\{\leq 1\}X_2\{\leq 1\}X_3\{\leq 1\}X_4\{\leq 0\}$ | $2*2*2*1 = 8$ | $\binom{4}{1} = 4$ |
| 0 | 1 | $[1, 0, 0, 0]$ | $X_1\{\leq 1\}X_2\{\leq 0\}X_3\{\leq 0\}X_4\{\leq 0\}$ | $2*1*1*1 = 2$ | $\binom{4}{1} = 4$ |

| $X_1$ | | 0 | | | 1 | | | 2 | | | | |
|---|---|---|---|---|---|---|---|---|---|---|---|---|
| $X_2$ | 0 | 1 | 2 | 0 | 1 | 2 | 0 | 1 | 2 | | | |
| | | | | | | | | | | 0 | 0 | |
| | | | | | | | | | | 1 | 1 | |
| | | | | | | | | | 1 | 2 | | |
| | | | | | | | | | | 0 | | |
| | | | | | | | | | 1 | 1 | 1 | |
| | | | | | | 1 | | 1 | 1 | 2 | | |
| | | | | | | | | | 1 | 0 | | |
| | | | | | | 1 | | 1 | 1 | 1 | 2 | |
| | | | 1 | | 1 | 1 | 1 | 1 | 1 | 2 | | |
| | | | | | | | | | | $X_4$ | $X_3$ | |

$$S\{\geq 3\} = S\{3\}$$

Fig. 6a. Binary success at level 3 with MUV cells highlighted: (a) in yellow for (2, 2, 2, 0)-type, and (b) in green for (2, 2, 1, 1)-type.

| $X_1$ | 0 | | | 1 | | | 2 | | | | |
|---|---|---|---|---|---|---|---|---|---|---|---|
| $X_2$ | 0 | 1 | 2 | 0 | 1 | 2 | 0 | 1 | 2 | | |
| | | | | | | | | | | 0 | |
| | | | | | | | | | | 1 | 0 |
| | | | | | | | | | 1 | 2 | |
| | | | | | | | | | | 0 | |
| | | | | | | | | | 1 | 1 | 1 |
| | | | | | 1 | | | 1 | 1 | 2 | |
| | | | | | | | | | 1 | 0 | |
| | | | | | 1 | | | 1 | 1 | 1 | 2 |
| | | | 1 | | 1 | 1 | | 1 | 1 | 2 | |
| | | | | | | | | | | $X_4$ | $X_3$ |

There are four 3-cell loops of the form
$X_1\{\geq 0\} X_2\{\geq 2\} X_3\{\geq 2\} X_4\{\geq 2\} =$
$X_1\{0, 1, 2\} X_2\{2\} X_3\{2\} X_4\{2\}$

$S\{\geq 3\} = S\{3\}$

There are six 4-cell loops of the form
$X_1\{\geq 1\} X_2\{\geq 2\} X_3\{\geq 1\} X_4\{\geq 2\} =$
$X_1\{1, 2\} X_2\{2\} X_3\{1, 2\} X_4\{2\}$

Fig. 6b. Binary success at level 3 with two sample prime-implicant (PI) loops indicated: (a) a PI loop extending from the yellow (2, 2, 2, 0)-type MUV cell (0,2,2,2) to the all-2 cell, and (b) a PI loop extending from the green (2, 2, 1, 1)-type MUV cell (1,2,1,2) to the all-2 cell. The total number of PI loops is 4 + 6 = 10.

| $X_1$ | 0 | | | 1 | | | 2 | | | | | |
|---|---|---|---|---|---|---|---|---|---|---|---|---|
| $X_2$ | 0 | 1 | 2 | 0 | 1 | 2 | 0 | 1 | 2 | | | |
| | | | | | | | | | | | 0 | |
| | | | | | | | | | | | 1 | 0 |
| | | | | | | | | | 1 | | 2 | |
| | | | | | | | | | | | 0 | |
| | | | | | | | | | 1 | | 1 | 1 |
| | | | | | 1 | | | 1 | 1 | | 2 | |
| | | | | | | | | | 1 | | 0 | |
| | | | | | 1 | | | 1 | 1 | | 1 | 2 |
| | | | 1 | | 1 | 1 | 1 | 1 | 1 | | 2 | |
| | | | | | | | | | | | $X_4$ | $X_3$ |

$$S_{PRE}\{\geq 3\} = S_{PRE}\{3\}$$

$= X_1\{2\} X_2\{2\} X_3\{2\} X_4\{0, 1, 2\} \vee X_1\{2\} X_2\{2\} X_3\{0, 1\} X_4\{2\}$
$\quad \vee X_1\{2\} X_2\{0, 1\} X_3\{2\} X_4\{2\} \vee X_1\{0, 1\} X_2\{2\} X_3\{2\} X_4\{2\}$
$\quad \vee X_1\{2\} X_2\{2\} X_3\{1\} X_4\{1\} \vee X_1\{2\} X_2\{1\} X_3\{2\} X_4\{1\}$
$\quad \vee X_1\{1\} X_2\{2\} X_3\{2\} X_4\{1\} \vee X_1\{2\} X_2\{1\} X_3\{1\} X_4\{2\}$
$\quad \vee X_1\{1\} X_2\{2\} X_3\{1\} X_4\{2\} \vee X_1\{1\} X_2\{1\} X_3\{2\} X_4\{2\}$

Fig. 6c. Binary success at level 3 covered with ten disjoint (non-overlapping) loops, comprising: (a) four loops (an intact 3-cell loop plus three curtailed 2-cell loops) that inherit the ones whose MUVs are (2, 2, 2, 0)-type, and (b) six loops (each a curtailed 1-cell loop containing just the MUV) that inherit the ones whose MUVs are (2, 2, 1, 1)-type. This threshold function is a shellable one, with each of the original PI loops retained as a single disjoint loop.

| $X_1$ | 0 | | | 1 | | | 2 | | | | | |
|---|---|---|---|---|---|---|---|---|---|---|---|---|
| $X_2$ | 0 | 1 | 2 | 0 | 1 | 2 | 0 | 1 | 2 | | | |
| | | | | | | | | | | 0 | | |
| | | | | | | | | | | 1 | 0 | |
| | | | | | | | | | 1 | 2 | | |
| | | | | | | | | | | 0 | | |
| | | | | | | | | | 1 | 1 | 1 | |
| | | | | | | 1 | | | 1 | 1 | 2 | |
| | | | | | | | | | 1 | 0 | | |
| | | | | | | 1 | | | 1 | 1 | 1 | 2 |
| | | | 1 | | 1 | 1 | 1 | 1 | 1 | 2 | | |
| | | | | | | | | | | $X_4$ | $X_3$ | |

$$S_{PRE}\{\geq 3\} = S_{PRE}\{3\}$$

$$= X_1\{2\}\, X_2\{1,2\}\, X_3\{2\}\, X_4\{1,2\} \vee X_1\{2\}\, X_2\{1,2\}\, X_3\{1\}\, X_4\{2\}$$
$$\vee X_1\{1\}\, X_2\{1,2\}\, X_3\{2\}\, X_4\{2\} \vee X_1\{2\}\, X_2\{2\}\, X_3\{1\}\, X_4\{1\}$$
$$\vee X_1\{1\}\, X_2\{2\}\, X_3\{1\}\, X_4\{2\} \vee X_1\{1\}\, X_2\{2\}\, X_3\{2\}\, X_4\{1\}$$
$$\vee X_1\{2\}\, X_2\{2\}\, X_3\{2\}\, X_4\{0\} \vee X_1\{2\}\, X_2\{2\}\, X_3\{0\}\, X_4\{2\}$$
$$\vee X_1\{2\}\, X_2\{0\}\, X_3\{2\}\, X_4\{2\} \vee X_1\{0\}\, X_2\{2\}\, X_3\{2\}\, X_4\{2\}$$

Fig. 6d. Binary success at level 3 with another coverage with alternative ten disjoint (non-overlapping) loops, comprising: (a) six loops (an intact 4-cell loop plus two curtailed 2-cell loops and three curtailed 1-cell loops) that inherit the ones whose MUVs are (2, 2, 1, 1)-type, and (b) four loops (each a curtailed 1-cell loop containing just the MUV) that inherit the ones whose MUVs are (2, 2, 2, 0)-type. This confirms the earlier observation that this threshold function is a shellable one, with each of the original PI loops retained as a single disjoint loop.

| $X_1$ | 0 | | | 1 | | | 2 | | | | |
|---|---|---|---|---|---|---|---|---|---|---|---|
| $X_2$ | 0 | 1 | 2 | 0 | 1 | 2 | 0 | 1 | 2 | | |
| | | | | | | | | | 1 | 0 | |
| | | | | | | 1 | | 1 | 1 | 1 | 0 |
| | | | 1 | | 1 | 1 | 1 | 1 | 1 | 2 | |
| | | | | | | 1 | | 1 | 1 | 0 | |
| | | | 1 | | 1 | 1 | 1 | 1 | 1 | 1 | 1 |
| | | 1 | 1 | 1 | 1 | 1 | 1 | 1 | 1 | 2 | |
| | | | 1 | | 1 | 1 | 1 | 1 | 1 | 0 | |
| | | 1 | 1 | 1 | 1 | 1 | 1 | 1 | 1 | 1 | 2 |
| | 1 | 1 | 1 | 1 | 1 | 1 | 1 | 1 | 1 | 2 | |
| | | | | | | | | | | $X_4$ | $X_3$ |

$$S\{\geq 2\} = S\{2,3\}$$

Fig.7a. Binary success at level 2 with MUV cells highlighted: (a) in yellow for (2, 2, 0, 0)-type, (b) in green for (2, 1, 1, 0)-type, and (c) in purple for (1, 1, 1, 1)-type.

| $X_1$ | 0 | | | 1 | | | 2 | | | | | |
|---|---|---|---|---|---|---|---|---|---|---|---|---|
| $X_2$ | 0 | 1 | 2 | 0 | 1 | 2 | 0 | 1 | 2 | | | |
| | | | | | | | | | 1 | 0 | | |
| | | | | | | 1 | | 1 | 1 | 1 | 0 | |
| | | | 1 | | 1 | 1 | 1 | 1 | 1 | 2 | | |
| | | | | | | 1 | | 1 | 1 | 0 | | |
| | | | 1 | | 1 | 1 | 1 | 1 | 1 | 1 | 1 | |
| | | 1 | 1 | 1 | 1 | 1 | 1 | 1 | 1 | 2 | | |
| | | | 1 | | 1 | 1 | 1 | 1 | 1 | 0 | | |
| | | 1 | 1 | 1 | 1 | 1 | 1 | 1 | 1 | 1 | 2 | |
| | 1 | 1 | 1 | 1 | 1 | 1 | 1 | 1 | 1 | 2 | | |
| | | | | | | | | | | $X_4$ | $X_3$ | |

$$S\{\geq 2\} = S\{2, 3\}$$

Fig. 7b. Binary success at level 2 with three sample prime-implicant (PI) loops indicated: (a) a 9-cell red PI loop extending from the yellow (2, 2, 0, 0)-type MUV cell (2,0,2,0) to the all-2 cell, (b) a shaded 12-cell PI loop extending from the (2, 1, 1, 0)-type MUV cell (0,1,1,2) to the all-2 cell, and (c) a green 16-cell PI loop extending from the purple (1 1, 1, 1)-type MUV cell (1,1,1,1) to the all-2 cell. The total number of PI loops is 6 + 12 + 1 = 19.

| $X_1$ | 0 | | | 1 | | | 2 | | | | | |
|---|---|---|---|---|---|---|---|---|---|---|---|---|
| $X_2$ | 0 | 1 | 2 | 0 | 1 | 2 | 0 | 1 | 2 | | | |
| | | | | | | | | | 1 | 0 | | |
| | | | | | | 1 | | 1 | 1 | 1 | 0 | |
| | | | 1 | | 1 | 1 | 1 | 1 | 1 | 2 | | |
| | | | | | | 1 | | 1 | 1 | 0 | | |
| | | | 1 | | 1 | 1 | 1 | 1 | 1 | 1 | 1 | |
| | | 1 | 1 | 1 | 1 | 1 | 1 | 1 | 1 | 2 | | |
| | | | 1 | | 1 | 1 | 1 | 1 | 1 | 0 | | |
| | | 1 | 1 | 1 | 1 | 1 | 1 | 1 | 1 | 1 | 2 | |
| | 1 | 1 | 1 | 1 | 1 | 1 | 1 | 1 | 1 | 2 | | |
| | | | | | | | | | | $X_4$ | $X_3$ | |

$$S_{PRE}\{\geq 2\} = S_{PRE}\{2, 3\}$$

$= X_1\{1, 2\} X_2\{1, 2\} X_3\{1, 2\} X_4\{1, 2\} \vee X_1\{0\} X_2\{1, 2\} X_3\{1, 2\} X_4\{2\}$
$\quad \vee X_1\{2\} X_2\{1, 2\} X_3\{0\} X_4\{1, 2\} \vee X_1\{1, 2\} X_2\{0\} X_3\{2\} X_4\{1, 2\}$
$\quad \vee X_1\{2\} X_2\{1, 2\} X_3\{1, 2\} X_4\{0\} \vee X_1\{0\} X_2\{1, 2\} X_3\{2\} X_4\{1\}$
$\quad \vee X_1\{1\} X_2\{1, 2\} X_3\{2\} X_4\{0\} \vee X_1\{2\} X_2\{0\} X_3\{1\} X_4\{1, 2\}$
$\quad \vee X_1\{1\} X_2\{2\} X_3\{0\} X_4\{1, 2\} \vee X_1\{0\} X_2\{2\} X_3\{1\} X_4\{1\}$
$\quad \vee X_1\{1\} X_2\{0\} X_3\{1\} X_4\{2\} \vee X_1\{1\} X_2\{1\} X_3\{0\} X_4\{2\}$
$\quad \vee X_1\{1\} X_2\{2\} X_3\{1\} X_4\{0\} \vee X_1\{0\} X_2\{0\} X_3\{2\} X_4\{2\}$
$\quad \vee X_1\{0\} X_2\{2\} X_3\{0\} X_4\{2\} \vee X_1\{0\} X_2\{2\} X_3\{2\} X_4\{0\}$
$\quad \vee X_1\{2\} X_2\{0\} X_3\{0\} X_4\{2\} \vee X_1\{2\} X_2\{0\} X_3\{2\} X_4\{0\}$
$\quad \vee X_1\{2\} X_2\{2\} X_3\{0\} X_4\{0\}$

Fig. 7c. Binary success at level 2 covered with 19 disjoint (non-overlapping) loops, comprising: (a) a single intact 16-cell loop, whose MUV is the (1, 1, 1, 1)-type, (b) 12 curtailed loops (four 4-cell loops, four 2-cell loops, and four 1-cell loops) that inherit the ones whose MUVs are (2, 1, 1, 0)-type, and (c) six loops (each a curtailed 1-cell loop containing just the MUV) that inherit the ones whose MUVs are (2, 2, 0, 0)-type. This threshold function is a shellable one, with each of the original PI loops retained as a single disjoint loop.

| $X_1$ | 0 | | | 1 | | | 2 | | | | |
|---|---|---|---|---|---|---|---|---|---|---|---|
| $X_2$ | 0 | 1 | 2 | 0 | 1 | 2 | 0 | 1 | 2 | | |
| | | | 1 | | 1 | 1 | 1 | 1 | 1 | 0 | |
| | | 1 | 1 | 1 | 1 | 1 | 1 | 1 | 1 | 1 | 0 |
| | 1 | 1 | 1 | 1 | 1 | 1 | 1 | 1 | 1 | 2 | |
| | | 1 | 1 | 1 | 1 | 1 | 1 | 1 | 1 | 0 | |
| | 1 | 1 | 1 | 1 | 1 | 1 | 1 | 1 | 1 | 1 | 1 |
| | 1 | 1 | 1 | 1 | 1 | 1 | 1 | 1 | 1 | 2 | |
| | 1 | 1 | 1 | 1 | 1 | 1 | 1 | 1 | 1 | 0 | |
| | 1 | 1 | 1 | 1 | 1 | 1 | 1 | 1 | 1 | 1 | 2 |
| | 1 | 1 | 1 | 1 | 1 | 1 | 1 | 1 | 1 | 2 | |
| | | | | | | | | | | $X_4$ | $X_3$ |

$$S\{\geq 1\} = S\{1, 2, 3\}$$

Fig. 8a. Binary success at level 1 with MUV cells highlighted: (a) in yellow for (2, 0, 0, 0)-type, and (b) in green for (1, 1, 0, 0)-type.

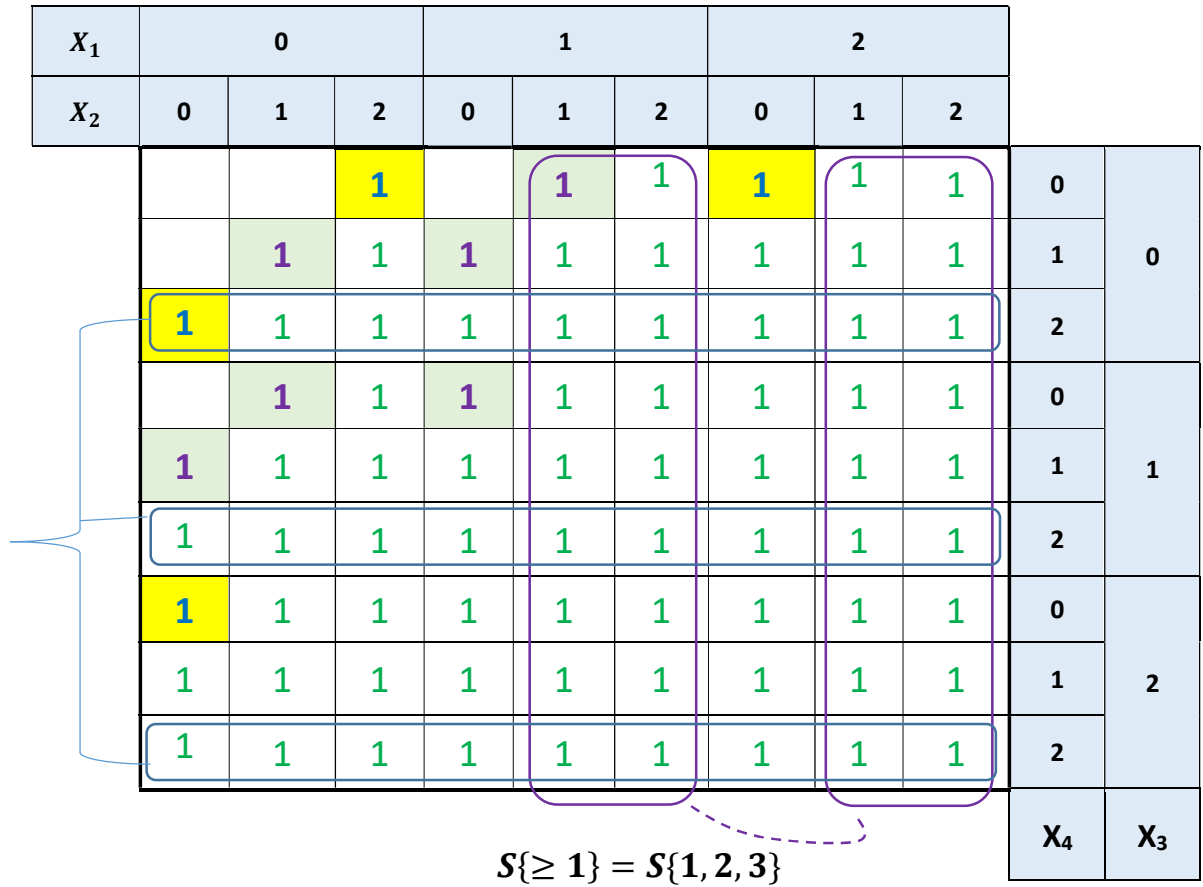

Fig. 8b. Binary success at level 1 with two sample prime-implicant (PI) loops indicated: (a) a 27-cell PI loop extending over 3 rows from the yellow (2, 0, 0, 0)-type MUV cell (0,0,0,2) to the all-2 cell, and (b) a 36-cell PI loop extending over 4 columns from the green (1, 1, 0, 0)-type MUV cell (1,1,0, 0) to the all-2 cell. The total number of PI loops is 4 + 6 = 10.

| $X_1$ | 0 | | | 1 | | | 2 | | | | |
|---|---|---|---|---|---|---|---|---|---|---|---|
| $X_2$ | 0 | 1 | 2 | 0 | 1 | 2 | 0 | 1 | 2 | | |
| | | | 1 | | 1 | 1 | 1 | 1 | 1 | 0 | |
| | | 1 | 1 | 1 | 1 | 1 | 1 | 1 | 1 | 1 | 0 |
| | 1 | 1 | 1 | 1 | 1 | 1 | 1 | 1 | 1 | 2 | |
| | | 1 | 1 | 1 | 1 | 1 | 1 | 1 | 1 | 0 | |
| | 1 | 1 | 1 | 1 | 1 | 1 | 1 | 1 | 1 | 1 | 1 |
| | 1 | 1 | 1 | 1 | 1 | 1 | 1 | 1 | 1 | 2 | |
| | 1 | 1 | 1 | 1 | 1 | 1 | 1 | 1 | 1 | 0 | |
| | 1 | 1 | 1 | 1 | 1 | 1 | 1 | 1 | 1 | 1 | 2 |
| | 1 | 1 | 1 | 1 | 1 | 1 | 1 | 1 | 1 | 2 | |
| | | | | | | | | | | $X_4$ | $X_3$ |

$$S_{PRE}\{\geq 1\} = S_{PRE}\{1, 2, 3\}$$
$$= X_1\{1, 2\}\, X_2\{0, 1, 2\}\, X_3\{1, 2\}\, X_4\{0, 1, 2\}$$
$$\vee\; X_1\{0\}\, X_2\{1, 2\}\, X_3\{1, 2\}\, X_4\{0, 1, 2\}$$
$$\vee\; X_1\{1, 2\}\, X_2\{0, 1, 2\}\, X_3\{0\}\, X_4\{1, 2\}$$
$$\vee\; X_1\{0\}\, X_2\{1, 2\}\, X_3\{0\}\, X_4\{1, 2\}$$
$$\vee\; X_1\{1, 2\}\, X_2\{1, 2\}\, X_3\{0\}\, X_4\{0\}$$
$$\vee\; X_1\{0\}\, X_2\{0\}\, X_3\{1, 2\}\, X_4\{1, 2\} \;\vee\; X_1\{0\}\, X_2\{0\}\, X_3\{0\}\, X_4\{2\}$$
$$\vee\; X_1\{0\}\, X_2\{0\}\, X_3\{2\}\, X_4\{0\} \;\vee\; X_1\{0\}\, X_2\{2\}\, X_3\{0\}\, X_4\{0\}$$
$$\vee\; X_1\{2\}\, X_2\{0\}\, X_3\{0\}\, X_4\{0\}$$

Fig. 8c. Binary success at level 1 covered with ten disjoint (non-overlapping) loops, comprising: (a) six loops (an intact 36-cell loop and two curtailed 12-cell loops plus three curtailed 4-cell loops) that inherit the ones whose MUVs are (1, 1, 0, 0)-type, and (b) four loops (each a curtailed 1-cell loop containing just the yellow MUV) that inherit the ones whose MUVs are (2, 0, 0, 0)-type. This threshold function is a shellable one, with each of the original PI loops retained as a single disjoint loop.

| $X_1$ | 0 | | | 1 | | | 2 | | | | |
|---|---|---|---|---|---|---|---|---|---|---|---|
| $X_2$ | 0 | 1 | 2 | 0 | 1 | 2 | 0 | 1 | 2 | | |
| | 1 | 1 | 1 | 1 | 1 | 1 | 1 | 1 | 1 | 0 | |
| | 1 | 1 | 1 | 1 | 1 | 1 | 1 | 1 | 1 | 1 | 0 |
| | 1 | 1 | 1 | 1 | 1 | 1 | 1 | 1 | | 2 | |
| | 1 | 1 | 1 | 1 | 1 | 1 | 1 | 1 | 1 | 0 | |
| | 1 | 1 | 1 | 1 | 1 | 1 | 1 | 1 | | 1 | 1 |
| | 1 | 1 | 1 | 1 | 1 | | 1 | | | 2 | |
| | 1 | 1 | 1 | 1 | 1 | 1 | 1 | 1 | | 0 | |
| | 1 | 1 | 1 | 1 | 1 | | 1 | | | 1 | 2 |
| | 1 | 1 | | 1 | | | | | | 2 | |
| | | | | | | | | | | $X_4$ | $X_3$ |

$$S\{\leq 2\} = S\{< 3\} = S\{0, 1, 2\}$$

Fig. 9a. Binary failure at level 3 with MLV cells highlighted: (a) in green for (2, 2, 1, 0)-type, and (b) in yellow for (2, 1, 1, 1)-type.

| $X_1$ | 0 | | | 1 | | | 2 | | | | |
|---|---|---|---|---|---|---|---|---|---|---|---|
| $X_2$ | 0 | 1 | 2 | 0 | 1 | 2 | 0 | 1 | 2 | | |
| | 1 | 1 | 1 | 1 | 1 | 1 | 1 | 1 | 1 | 0 | |
| | 1 | 1 | 1 | 1 | 1 | 1 | 1 | 1 | 1 | 1 | 0 |
| | 1 | 1 | 1 | 1 | 1 | 1 | 1 | 1 | | 2 | |
| | 1 | 1 | 1 | 1 | 1 | 1 | 1 | 1 | 1 | 0 | |
| | 1 | 1 | 1 | 1 | 1 | 1 | 1 | 1 | | 1 | 1 |
| | 1 | 1 | 1 | 1 | 1 | | 1 | | | 2 | |
| | 1 | 1 | 1 | 1 | 1 | 1 | 1 | 1 | | 0 | |
| | 1 | 1 | 1 | 1 | 1 | | 1 | | | 1 | 2 |
| | 1 | 1 | | 1 | | | | | | 2 | |
| | | | | | | | | | | $X_4$ | $X_3$ |

$$S_{PRE}\{\leq 2\} = S_{PRE}\{< 3\} = S_{PRE}\{0, 1, 2\}$$
$$= X_1\{0,1\}\, X_2\{0,1,2\}\, X_3\{0,1\}\, X_4\{0,1\}$$
$$\lor X_1\{2\}\, X_2\{0,1\}\, X_3\{0,1\}\, X_4\{0,1\} \lor X_1\{0,1\}\, X_2\{0,1\}\, X_3\{2\}\, X_4\{0,1\}$$
$$\lor X_1\{0,1\}\, X_2\{0,1\}\, X_3\{0,1\}\, X_4\{2\} \lor X_1\{0\}\, X_2\{0,1\}\, X_3\{2\}\, X_4\{2\}$$
$$\lor X_1\{2\}\, X_2\{2\}\, X_3\{0\}\, X_4\{0,1\} \lor X_1\{0\}\, X_2\{2\}\, X_3\{2\}\, X_4\{0,1\}$$
$$\lor X_1\{0\}\, X_2\{2\}\, X_3\{0,1\}\, X_4\{2\} \lor X_1\{2\}\, X_2\{0,1\}\, X_3\{0\}\, X_4\{2\}$$
$$\lor X_1\{2\}\, X_2\{0,1\}\, X_3\{2\}\, X_4\{0\} \lor X_1\{1\}\, X_2\{0\}\, X_3\{2\}\, X_4\{2\}$$
$$\lor X_1\{1\}\, X_2\{2\}\, X_3\{0\}\, X_4\{2\} \lor X_1\{1\}\, X_2\{2\}\, X_3\{2\}\, X_4\{0\}$$
$$\lor X_1\{2\}\, X_2\{0\}\, X_3\{1\}\, X_4\{2\} \lor X_1\{2\}\, X_2\{0\}\, X_3\{2\}\, X_4\{1\}$$
$$\lor X_1\{2\}\, X_2\{2\}\, X_3\{1\}\, X_4\{0\}$$

Fig. 9b. Binary failure at level 3 covered with 16 disjoint (non-overlapping) loops, comprising: (a) four loops (an intact 24-cell loop plus three curtailed 8-cell loops) that inherit the ones whose MLVs are (2, 1, 1, 1)-type, and (b) 12 loops (6 curtailed 2-cell loops plus 6 curtailed 1-cell loops (just the MLVs)) that inherit the ones whose MLVs are (2, 2, 1, 0)-type. This threshold function is a shellable one, with each of the original PI loops retained as a single disjoint loop.

| $X_1$ | 0 | | | 1 | | | 2 | | | | |
|---|---|---|---|---|---|---|---|---|---|---|---|
| $X_2$ | 0 | 1 | 2 | 0 | 1 | 2 | 0 | 1 | 2 | | |
| | 1 | 1 | 1 | 1 | 1 | 1 | 1 | 1 | | 0 | |
| | 1 | 1 | 1 | 1 | 1 | | 1 | | | 1 | 0 |
| | 1 | 1 | | 1 | | | | | | 2 | |
| | 1 | 1 | 1 | 1 | 1 | | 1 | | | 0 | |
| | 1 | 1 | | 1 | | | | | | 1 | 1 |
| | 1 | | | | | | | | | 2 | |
| | 1 | 1 | | 1 | | | | | | 0 | |
| | 1 | | | | | | | | | 1 | 2 |
| | | | | | | | | | | 2 | |
| | | | | | | | | | | $X_4$ | $X_3$ |

$$S\{\leq 1\} = S\{< 2\} = S\{0, 1\}$$

Fig. 10a. Binary failure at level 2 with MLV cells highlighted: (a) in green for (2, 1, 0, 0)-type, and (b) in yellow for (1, 1, 1, 0)-type.

| $X_1$ | 0 | | | 1 | | | 2 | | | | |
|---|---|---|---|---|---|---|---|---|---|---|---|
| $X_2$ | 0 | 1 | 2 | 0 | 1 | 2 | 0 | 1 | 2 | | |
| | 1 | 1 | 1 | 1 | 1 | 1 | 1 | 1 | | 0 | |
| | 1 | 1 | 1 | 1 | 1 | | 1 | | | 1 | 0 |
| | 1 | 1 | | 1 | | | | | | 2 | |
| | 1 | 1 | 1 | 1 | 1 | | 1 | | | 0 | |
| | 1 | 1 | | 1 | | | | | | 1 | 1 |
| | 1 | | | | | | | | | 2 | |
| | 1 | 1 | | 1 | | | | | | 0 | |
| | 1 | | | | | | | | | 1 | 2 |
| | | | | | | | | | | 2 | |
| | | | | | | | | | | $X_4$ | $X_3$ |

$$S_{PRE}\{\leq 1\} = S_{PRE}\{< 2\} = S_{PRE}\{0, 1\}$$
$$= X_1\{0\} X_2\{0, 1\} X_3\{0, 1\} X_4\{0, 1\} \vee X_1\{1\} X_2\{0, 1\} X_3\{0\} X_4\{0, 1\}$$
$$\vee X_1\{1\} X_2\{0, 1\} X_3\{1\} X_4\{0\} \vee X_1\{1\} X_2\{0\} X_3\{1\} X_4\{1\}$$
$$\vee X_1\{0\} X_2\{0, 1\} X_3\{0\} X_4\{2\} \vee X_1\{0\} X_2\{0, 1\} X_3\{2\} X_4\{0\}$$
$$\vee X_1\{0\} X_2\{2\} X_3\{0\} X_4\{0, 1\} \vee X_1\{2\} X_2\{0\} X_3\{0\} X_4\{0, 1\}$$
$$\vee X_1\{0\} X_2\{0\} X_3\{1\} X_4\{2\} \vee X_1\{0\} X_2\{0\} X_3\{2\} X_4\{1\}$$
$$\vee X_1\{0\} X_2\{2\} X_3\{1\} X_4\{0\} \vee X_1\{1\} X_2\{0\} X_3\{0\} X_4\{2\}$$
$$\vee X_1\{1\} X_2\{0\} X_3\{2\} X_4\{0\} \vee X_1\{1\} X_2\{2\} X_3\{0\} X_4\{0\}$$
$$\vee X_1\{2\} X_2\{0\} X_3\{1\} X_4\{0\} \vee X_1\{2\} X_2\{1\} X_3\{0\} X_4\{0\}$$

Fig. 10b. Binary failure at level 2 covered with 16 disjoint (non-overlapping) loops, comprising: (a) four loops (an intact 8-cell loop plus three curtailed 4-cell, 2-cell and 1-cell loops) that inherit the ones whose MLVs are (1, 1, 1, 0)-type, and (b) 12 loops (4 curtailed 2-cell loops plus 8 curtailed 1-cell loops (just the MLVs)) that inherit the ones whose MLVs are (2, 1, 0, 0)-type. This threshold function is a shellable one, with each of the original PI loops retained as a single disjoint loop.

| $X_1$ | 0 | | | 1 | | | 2 | | | | |
|---|---|---|---|---|---|---|---|---|---|---|---|
| $X_2$ | 0 | 1 | 2 | 0 | 1 | 2 | 0 | 1 | 2 | | |
| | 1 | 1 | | 1 | | | | | | 0 | |
| | 1 | | | | | | | | | 1 | 0 |
| | | | | | | | | | | 2 | |
| | 1 | | | | | | | | | 0 | |
| | | | | | | | | | | 1 | 1 |
| | | | | | | | | | | 2 | |
| | | | | | | | | | | 0 | |
| | | | | | | | | | | 1 | 2 |
| | | | | | | | | | | 2 | |
| | | | | | | | | | | $X_4$ | $X_3$ |

$$S\{\leq 0\} = S\{< 1\} = S\{0\}$$

Fig. 11a. Binary failure at level 1 with MLV cells highlighted in yellow for (1, 0, 0, 0)-type.

| $X_1$ | 0 | | | 1 | | | 2 | | | | |
|---|---|---|---|---|---|---|---|---|---|---|---|
| $X_2$ | 0 | 1 | 2 | 0 | 1 | 2 | 0 | 1 | 2 | | |
| | 1 | 1 | | 1 | | | | | | 0 | |
| | 1 | | | | | | | | | 1 | 0 |
| | | | | | | | | | | 2 | |
| | 1 | | | | | | | | | 0 | |
| | | | | | | | | | | 1 | 1 |
| | | | | | | | | | | 2 | |
| | | | | | | | | | | 0 | |
| | | | | | | | | | | 1 | 2 |
| | | | | | | | | | | 2 | |
| | | | | | | | | | | $X_4$ | $X_3$ |

$$S_{PRE}\{\leq 0\} = S_{PRE}\{< 1\} = S_{PRE}\{0\}$$
$$= X_1\{0\}\, X_2\{0, 1\}\, X_3\{0\}\, X_4\{0\} \vee X_1\{0\}\, X_2\{0\}\, X_3\{0\}\, X_4\{1\}$$
$$\vee X_1\{0\}\, X_2\{0\}\, X_3\{1\}\, X_4\{0\} \vee X_1\{1\}\, X_2\{0\}\, X_3\{0\}\, X_4\{0\}$$

Fig. 11b. Binary failure at level 1 covered with four disjoint (non-overlapping) loops, comprising an intact 2-cell loop plus three curtailed 1-cell loops that inherit the ones whose MLVs are (1, 0, 0, 0)-type.